\documentclass[10pt,a4paper]{article}

\usepackage{amsmath,amsthm}
\usepackage[authoryear]{natbib}
\usepackage{microtype}
\usepackage[bitstream-charter]{mathdesign}
\usepackage{bm}
\usepackage[final]{showkeys}
\usepackage{color}
\usepackage[title]{appendix}

\numberwithin{equation}{section}
\usepackage[hmargin=2.5cm,vmargin={2.5cm,3cm}]{geometry}
\allowdisplaybreaks

\newtheorem{theorem}{Theorem}[section]
\newtheorem{proposition}{Proposition}[section]
\theoremstyle{remark}
\newtheorem{remark}{Remark}[section]
\theoremstyle{definition}
\newtheorem{definition}{Definition}[section]

\title{Fractional Diffusion-Telegraph Equations and their Associated Stochastic Solutions}

\author{$\text{Mirko D'Ovidio}_1$, $\text{Federico Polito}_2$ \\
	\footnotesize (1) -- Dipartimento di Scienze di Base e Applicate per l'Ingegneria,
		``Sapienza'' Universit\`{a} di Roma\\
	\footnotesize Via A. Scarpa 16, 00161 Roma, Italy\\
	\footnotesize Email address: mirko.dovidio@uniroma1.it\\
	\footnotesize (2) -- Dipartimento di Matematica ``G.\ Peano'', Universit\`{a} degli Studi di Torino\\
	\footnotesize Via Carlo Alberto 10, 10123 Torino, Italy\\
	\footnotesize Email address: federico.polito@unito.it\\
	}
	
\date{}

\begin{document}

	\maketitle
	
	\begin{abstract}

		\noindent We present the stochastic solution to a generalized fractional partial differential equation
        involving a regularized operator related to the so-called Prabhakar operator and admitting, amongst others, as specific cases
        the fractional diffusion equation and the fractional telegraph equation. The stochastic solution
        is expressed as a L\'evy process time-changed with the inverse process to a linear combination
        of (possibly subordinated) independent stable subordinators of different indices.
        Furthermore a related SDE is derived and discussed.
        
        \medskip
        

	\end{abstract}
	
	\section{Introduction}
	
        In the last few decades considerable effort has been devoted
        to the study of fractional partial differential equations (fPDEs)
        that is of PDEs in which usual differential operators
        are substituted by fractional differential operators
        (for example two rather recent references are \citet{kochubei} and \citet{leonenko}).
        The simplest equation of this class is the so-called
        fractional diffusion equation, also known as diffusion-wave
        equation \citep[see amongst others][]{schneider,mainardi,meerschaert,ors1}.
        Another well-known and well-studied fPDE is the
        fractional telegraph equation \citep{ors2,camargo,povstenko,yakubovich,DTO2012}.
        In the more recent years, moreover, an increasing number of papers presented
        results connecting the study of fractional PDEs to
        that of some time-changed stochastic processes. The aim of this paper is
        to clarify this connection for a very general class of fPDEs which includes
        as specific cases both parabolic and hyperbolic fPDEs as well as
        more general integral and differential equations.
        In order to be more specific and for the sake of comprehension, we will start by recalling here
        the definitions of the classical fractional operators of Riemann--Liouville type
        and the Dzhrbashyan--Caputo derivative (see for the latter \citet{liouville,nersesyan,caputo} --- see \citet{kilbas} for a reference book).
        
   		\begin{definition}[Riemann--Liouville integral]
		    Let $f \in L^1_{\text{loc}}(0,b)$, $0 < t < b \le \infty$, be a locally integrable real-valued function.
		    The operator
		    \begin{align}
		        \label{rlint}
		        J^{\alpha}_t f & =\frac{1}{\Gamma(\alpha)}\int_{0}^t\frac{f(u)}{(t-u)^{1-\alpha}}\textup{d}u, \qquad \alpha > 0,
		    \end{align}
		    is called Riemann--Liouville integral of order $\alpha$.
		\end{definition}

		\begin{definition}[Riemann--Liouville derivative]
			\label{rie}
		    Let $f \in L^1(0,b)$, $-\infty \le a < t < b \le \infty$, and
		    define the power-law kernel $\mathcal{L}_\beta (t) = t^{\beta-1}/\Gamma(\beta)$, $\beta>0$.
		    Consider $\alpha > 0$ and write $m = \lceil \alpha \rceil$ for the smallest integer greater than or equal to
		    $\alpha$. For $f \ast \mathcal{L}_{m-\alpha} \in W^{m,1}(0,b)$,
		    where $W^{m,1}(0,b)$ is the Sobolev space
		    \begin{align}
		        W^{m,1}(0,b) = \left\{ h \in L^1(0,b) \colon \frac{\textup{d}^m}{\textup{d}t^m} h \in L^1(0,b) \right\},
		    \end{align}
		    the Riemann--Liouville derivative of order $\alpha$ is defined as
		    \begin{align}
		        \label{rlder}
		        \frac{\textup{d}^\alpha}{\textup{d} t^\alpha}f(t) = \frac{1}{%
		        \Gamma(m-\alpha)} \frac{\textup{d}^m}{\textup{d}t^m} \int_{0}^t (t-s)^{m-1-\alpha}f(s) \, \textup{d}s.
		    \end{align}
		\end{definition}

		In order to introduce the definition of the Dzhrbashyan--Caputo derivative, let us
		denote by $AC^{n}\left(0,b\right)$, $n\in \mathbb{N}$, the
		space of real-valued functions $h\left( t\right)$ with
		continuous derivatives up to order $n-1$ on $\left(0,b\right)$
		and such that $h^{\left(n-1\right) }\left(t\right)$ belongs to the space of absolutely continuous functions
		$AC\left(0,b\right)$, i.e.\
		\begin{equation}
			AC^{n}\left(0,b\right) =\left\{ h : \left(0,b\right) \rightarrow
			\mathbb{R}\colon\frac{\textup{d}^{n-1}}{\textup{d}x^{n-1}}f \left( x\right) \in AC\left(
			0,b\right) \right\} .
		\end{equation}

		\begin{definition}[Dzhrbashyan--Caputo derivative]
			\label{capu}
		    Let $\alpha>0$, $m = \lceil \alpha \rceil$. 
		    The Dzhrbashyan--Caputo derivative of order $\alpha>0$ is defined as
		    \begin{equation}
		        \label{Capu}
		        \frac{\mathfrak{d}^\alpha}{\mathfrak{d}t^\alpha} f(t)= \frac{1%
		        }{\Gamma(m-\alpha)}\int_a^{t}(t-s)^{m-1-\alpha}\frac{\textup{d}^m}{\textup{d}s^m}f(s) \, \textup{d}s
		    \end{equation}
		    for $f \in AC^m(0,b)$ such that \eqref{Capu} exists.
		\end{definition}       
        
        Let us now recall one of the most famous applications of the Dzhrbashyan--Caputo derivative:
        the fractional diffusion equation in dimension one in its simplest form.
        Let us consider thus the Cauchy problem
        \begin{align}
            \label{meme}
            \begin{cases}
                \frac{\mathfrak d^\alpha}{\mathfrak{d} t^\alpha} u(x,t) = \lambda^2 \frac{\mathrm d^2}{\mathrm dx^2} u(x,t), & t > 0,
                \: x \in \mathbb{R}, \\
                u(x,0) = \delta(x), & 0 < \alpha \le 2, \\
                \left. \frac{\mathrm d}{\mathrm dt} u(x,t) \right|_{t=0} = 0, & 1 < \alpha \le 2.
            \end{cases}
        \end{align}
        It has been proven that the solution to \eqref{meme} can be written as \citep{ors1,kilbas}
        \begin{align}
            \label{so}
            \frac{1}{2\lambda t^{\alpha/2}} W_{-\alpha/2, 1-\alpha/2} \left( -\frac{|x|}{\lambda t^{\alpha/2}} \right), \qquad
            t \ge 0, \: x \in \mathbb{R}, \: 0< \alpha \le 2,
        \end{align}
        where $W_{a,b}(z)$ is the Wright function \citep[Chapter 1]{kilbas}. The solution \eqref{so} has the remarkable
        property that it reduces to the Gaussian function for $\alpha=1$ and to the classical d'Alambert's solution to
        the wave equation for $\alpha \to 2$ while keeping an intermediate behaviour for $\alpha \in (1,2)$.
        
        Aside the analytical point of view, by starting from the well-known fact that diffusion processes
        are strictly connected via their distributional
        structure to parabolic equations, the most part of the recent research on the subject has been in fact dedicated to
        construct stochastic processes that can be related to various classes of fractional PDEs
        and that can therefore furnish them with a microscopic interpretation.
        For example the fractional diffusion equation \eqref{meme} can be related to a time-changed
        Brownian motion (see e.g.\ \citet{meerschaert,ors1,meerschaert2004limit,meerschaert2011stochastic}).
        Indeed, let us call $B_t$, $t \ge 0$, a standard Brownian motion and $V_t^\alpha$, $t \ge 0$, $\alpha \in (0,1)$ an
        $\alpha$-stable subordinator, independent of $B_t$, from which the time-change will be constructed.
        Note that stable subordinators are a particularly well-behaved class of L\'evy processes which are increasing
        and have a very simple Laplace exponent. Then, after having defined the
        right-inverse process as $K_t^\alpha = \inf \{s \ge 0 \colon V_s^\alpha > t \}$, $t \ge 0$
        (see \citet{bingham1971limit} for details) we have that the marginal distribution
        $\mathbb{P} \{ B_{K_t^\alpha} \in \mathrm dx \}/\mathrm dx$ is the solution to \eqref{meme}.
        Similar considerations can be done for hyperbolic PDEs and some stochastic
        processes describing particles moving with finite velocity. With respect to this the reader can consult the papers by
        \citet{ors2,kochu,DTO2012,gop} and the references therein. For what concerns fractional evolution equations
        in abstract spaces we refer to \citet{baeumer2001stochastic,meerschaert2009fractional,eidelman2004cauchy,
        bazhlekova2000subordination,bazhlekova}.
        
        In this paper we study fractional PDEs in which the operator acting on time generalizes
        the Dzhrbashyan--Caputo fractional derivative and connect them to time-changed L\'evy processes.
        The general fractional PDE that we study in Section \ref{lele} generalizes both diffusion-like and telegraph-like differential equations.
        
        Let us thus start by considering the operator $\bm{\mathrm D}_{\alpha,\eta,\zeta;0+}^\xi$ that we call
        Prabhakar derivative. It was first
		introduced by \citet{saigo} and it is defined as
		\begin{align}
			\label{00tex}
			\left( \bm{\mathrm D}_{\alpha,\eta,\zeta;0+}^\xi f(\cdot) \right) (t) = \frac{\mathrm d^{\eta+\theta}}{\mathrm dt^{\eta+\theta}}
			\int_0^t (t-y)^{\theta-1} E_{\alpha,\theta}^{-\xi}\left[ \zeta (t-y)^\alpha \right]
			f(y) \, \mathrm dy,
		\end{align}
		where $\theta > 0$, $\eta > 0$, $\zeta \in \mathbb{R}$, $t \ge 0$ and  
		\begin{align}
			\label{00poch}
			E_{\alpha,\eta}^\xi (x) = \sum_{r=0}^\infty \frac{x^r(\xi)_r}{r! \Gamma(\alpha r + \eta)},
			\qquad \alpha,\eta,\xi \in \mathbb{R}, \: \alpha> 0,
		\end{align}
		is known as the generalized Mittag--Leffler function (see for example \citet{saigo} or \citet{prab}).
		The symbol 
		$(\xi)_r$ in \eqref{00poch} is the so-called Pochhammer symbol.
		Recall that the operator in
		\eqref{00tex} is the Riemann--Liouville derivative and notice that the operator $\bm{\mathrm D}_{\alpha,\eta,\zeta;0+}^\xi$ is the
		left-inverse to the convolution-type operator \citep{saigo}
		\begin{align}
			\label{20ciao}
			\left( \bm{\mathrm E}_{\alpha,\eta,\zeta;0+}^\xi f(\cdot) \right) (t)
			= \int_0^t (t-y)^{\eta-1} E_{\alpha,\eta}^\xi \left[ \zeta(t-y)^\alpha \right] f(y) \, \mathrm dy,
		\end{align}
		which was originally introduced by \citet{prab}.
        Consider now the following operator:
		\begin{definition}[Regularized Prabhakar derivative]      
        	Let $\eta>0$, $\xi,\zeta \in \mathbb{R}$, $\alpha> 0$, $m=\lceil \eta \rceil$, $\kappa = \lceil \xi \rceil$,
        	$f \in \mathbb{A}^{\kappa} (\mathbb{R}^+)$, where
    		\begin{align}
    			\mathbb{A}^{\kappa}(\mathbb{R}^+) = \biggl\{ u \colon \mathbb{R}^+ \mapsto
                \mathbb{R}^+ \text{ s.t. } \sum_{j=0}^{\kappa-1} a_j\frac{\mathfrak{d}^{\eta - \alpha j}}{\mathfrak{d} t^{\eta - \alpha j}}
                u \in C(\mathbb{R}^+), \; \alpha \in (0,1], \: a_j >0 \; \forall\, j, \:
      		    \big| \frac{\textup{d} u}{\textup{d} t}(t) \big| \leq t^{\beta-1}, \; \beta>0	\biggr\}. \notag
	    	\end{align}
	    	The operator
			\begin{align}
				\label{vtvt}
				\left( \mathbb D_{\alpha,\eta,\zeta;0+}^\xi f(\cdot) \right) (t)
				= \left( \bm{\mathrm D}_{\alpha,\eta,\zeta;0+}^\xi f(\cdot) \right) (t)
				- \sum_{k=0}^{m-1} f^{(k)}(0^+) t^{k-\eta} E^{-\xi}_{\alpha,k-\eta+1}(\zeta t^{\alpha}),
			\end{align}
			is called regularized Prabhakar derivative.
	    \end{definition}	
		Note that the regularized Prabhakar derivative has Laplace transform
		\begin{align}
			\label{barrique}
			\int_0^\infty e^{-st} \left( \mathbb D_{\alpha,\eta,\zeta;0+}^\xi f(\cdot) \right) (t) \, \mathrm dt
			= s^\eta (1-\zeta s^{-\alpha})^\xi \tilde{f}(s) -
			\sum_{k=0}^{m-1}				
			f^{(k)}(0^+) \, s^{\eta-k-1}(1-\zeta s^{-\alpha})^\xi,
		\end{align}
		where $\tilde{f}(s)$ is the Laplace transform of function $f$.

        In Section \ref{george} we introduce a stochastic process $M_t^\delta$, $t \ge 0$, $\delta \ge 0$,
        built as a suitable linear combination of stable subordinators that are made dependent by a common random time-change with a further
        independent stable subordinator. Indeed formula \eqref{ddir} tells us that $M^\delta_t = 
			\sum_{r=0}^n \binom{n}{r}^{\delta/(n\omega_r)} {}_r V_{t}^{\omega_r}$, $t \ge 0$, where ${}_r V_t^{\omega_r}$,
		$r=1,\dots n$, are the dependent $\omega_r$-stable subordinators with $\omega_r = \gamma+\nu -r \nu \delta/n$,
		$n= \lceil \delta \rceil$ and $\nu \delta < \gamma+\nu < 1$.    
        The main result of Section \ref{george} is Theorem \ref{55teo}
        which states that the process $Z_t^\delta := \inf \left\{ s \ge 0 \colon M^\delta_t \notin (0,t) \right\}$, $t \ge 0$,
        i.e.\ the right-inverse process to $M_t^\delta$, has marginal distribution which solves the Cauchy problem
		\begin{align}
			\begin{cases}
				\left( \mathbb{D}_{\nu,\gamma+\nu,-1;0+}^\delta h(x,\cdot) \right) (t)
				= -\frac{\partial}{\partial x} h(x,t), \qquad t\ge 0, \: x \ge 0, \\
				h(x,0^+) = \delta(x).
			\end{cases}
		\end{align}       
        
        In Section \ref{lele} the most general fPDE we deal with is
   		\begin{align}
            \label{astratto}
			\left( \mathbb{D}_{\nu,\gamma+\nu,-1;0+}^\delta g(x,\cdot) \right) (t)
			= \mathcal{A}g(x, t), \qquad x \in \mathbb{R}^d, \: t>0,
		\end{align}
        with $\delta > 0$, $\delta \nu < \gamma+\nu \le 1$, $\gamma,\nu \in (0,1)$,
        and where $\mathcal{A}$ is the infinitesimal generator of a L\'evy process $A_t^x$, $t \ge 0$, starting
        from $x \in \mathbb{R}^d$. We prove in Theorem \ref{TC-telegraph} that the time-changed process
        $A^x_{Z_t^\delta}$, $t \ge 0$,
        has one-dimensional distribution which solves \eqref{astratto} with a suitable choice of the initial datum.
        
        Section \ref{last} contains all the details for the case $\mathcal{A} = \partial^2 / \partial x^2$.
        However, in that section $\delta \nu < \gamma+\nu \le 2$, hence including important specific cases such as
        the classical telegraph equation. This clearly does not permit us to relate the solution to $A^x_{Z_t^\delta}$
        but still several results from an analytical point of view are obtained.

    \section{The operator $\bm{\mathbb D_{\alpha,\eta,\zeta;0+}^\xi}$ and its relation to the hitting time of
		linear combinations of stable subordinators}
        
        \label{george}
		In the following we study the connection
		of the operator $\mathbb D_{\alpha,\eta,\zeta;0+}^\xi$ to some
		stochastic processes constructed as inverse processes to some linear combinations of
		dependent stable subordinators.				
				
		\subsection{Connections of $\bm{\mathbb D_{\alpha,\eta,\zeta;0+}^\xi}$ to the hitting time of
			linear combinations of stable subordinators}
		
			Here we analyze the relations of some inverse processes of linear combinations
			of dependent stable subordinators to differential equations involving the operator
			$\mathbb D_{\alpha,\eta,\zeta;0+}^\xi$ with a restriction on the range of the parameter $\eta$ and 
			considered as an operator acting on functions $t \mapsto f(\cdot, t)$. When $\eta \in (0,1)$ from formula \eqref{vtvt} we have
			\begin{align}
				\left( \mathbb{D}_{\alpha,\eta,\zeta;0+}^\xi f(x,\cdot) \right) (t)
				= \left( \bm{\mathrm D}_{\alpha,\eta,\zeta;0+}^\xi f(x,\cdot) \right) (t) - f(x, 0^+)
				\, t^{-\eta} E_{\alpha,1-\eta}^{-\xi} (\zeta t^\alpha),
			\end{align}
			with Laplace transform
			\begin{align}
				\label{LapCaputo-typeOp}
				\int_0^\infty e^{-s t} \left( \mathbb{D}_{\alpha,\eta,\zeta;0+}^\xi f(x,\cdot) \right)(t)\, \mathrm dt
				= s^{\eta}(1-\zeta s^{-\alpha})^{\xi}\tilde{f}(x, s) - f(x, 0^+)
				s^{\eta-1}(1-\zeta s^{-\alpha})^{\xi}.
			\end{align}
			
			Let us now assume $(\gamma+\nu) \in (0,1)$, $\delta \in \mathbb{R}$, $k=\lceil \delta\rceil$, and
			call $\Omega^{d+1} := \Omega \times(0,\infty)$ with $\Omega \subseteq \mathbb{R}^d$.
			We define the function space
    		\begin{align}
    			\label{rafaela}
    			\mathbb{A}^{k}(\Omega^{d+1}) = \biggl\{ u \colon \Omega^{d+1} \mapsto
                \mathbb{R}^+ \text{ s.t. } & \sum_{j=0}^{k-1} a_j\frac{\mathfrak{d}^{\eta - \alpha j}}{\mathfrak{d} t^{\eta - \alpha j}}
                u \in C(\Omega^{d+1}), \; \alpha \in (0,1], \; a_j >0 \; \forall\, j,\\
      		    & \big| \frac{\partial u}{\partial t}(x,t) \big| \leq \mathfrak{g}(x)t^{\beta-1}, \; \beta>0,\; \mathfrak{g} \in L^\infty(\Omega)    
	    		\biggr\}. \notag 
	    	\end{align}
			Our aim is to find the solution to the following Cauchy problem, $u \in \mathbb{A}^{k}(\mathbb{R}^{1+1}_+) \cap C_0(\mathbb{R}_+)$,
			\begin{align}
				\label{04sub}
				\begin{cases}
					\left( \mathbb{D}_{\nu,\gamma+\nu,-1;0+}^\delta h(x,\cdot) \right) (t)
					= -\frac{\partial}{\partial x} h(x,t), \qquad t > 0, \: x > 0, \\
					h(x,0^+) = \delta(x).
				\end{cases}
			\end{align}
			It is interesting to note that when $\delta=0$, the above equation reduces to that satisfied by the density law of the
			inverse $(\gamma+\nu)$-stable subordinator (see e.g.\ \citet[Section 4]{mirko2}). Indeed, by
			using Remark \ref{valmora} and after some calculations, we arrive at
			\begin{align}
				\frac{\mathfrak{d}^{\gamma+\nu}}{\mathfrak{d} t^{\gamma+\nu}} h(x,t) = -\frac{\partial}{\partial x}
				h(x,t), \qquad t > 0, \: x > 0,
			\end{align}
			where, we recall, $\mathfrak{d}^\varpi / \mathfrak{d} t^\varpi$
			is the Dzhrbashyan--Caputo fractional deivative of order $\varpi$.
			For $\delta=1$, we obtain
			\begin{align}
				\label{adv}
				\frac{\mathfrak{d}^{\gamma+\nu}}{\mathfrak{d} t^{\gamma+\nu}} h(x,t)
				+ \frac{\mathfrak{d}^\gamma}{\mathfrak{d} t^\gamma} h(x,t) = -\frac{\partial}{\partial x}
				h(x,t), \qquad t > 0, \: x > 0.
			\end{align}
	
			In order to solve \eqref{04sub} we apply a Laplace--Laplace transform to \eqref{04sub}
			with respect to both variables.
			\begin{proposition}
				The Laplace--Laplace transform $\tilde{\tilde{h}}(z,s)
				= \int_0^\infty \int_0^\infty e^{-zx-st} h(x,t) \, \mathrm dx \, \mathrm dt$
				of the solution to \eqref{04sub} is
				\begin{align}
					\label{hlaplap}
					\tilde{\tilde{h}}(z,s) = \frac{s^{\gamma+\nu-1}\left( 1+s^{-\nu} \right)^\delta}{
					s^{\gamma+\nu} \left( 1+s^{-\nu} \right)^\delta + z}, \qquad z>0, \: s>0.
				\end{align}
			
				\begin{proof}
					By means of direct calculation and considering \eqref{barrique}
					and the initial condition we have that
					\begin{align}
						\label{ebib}
						s^{\gamma+\nu} \tilde{\tilde{h}}(z,s) \left( 1+s^{-\nu} \right)^\delta
						= -z \, \tilde{\tilde{h}}(z,s) + s^{\gamma+\nu-1} \left( 1+s^{-\nu} \right)^\delta.
					\end{align}
					Rewriting \eqref{ebib} with respect to $\tilde{\tilde{h}}(z,s)$ we obtain
					the thesis \eqref{hlaplap}.
				\end{proof}
		\end{proposition}		
		
		\begin{proposition}
			The $t$-Laplace transform and $x$-Laplace transform of $h(x,t)$ read respectively
			\begin{align}
				& \label{05lap} \tilde{h}(x,s) = \int_0^\infty e^{-st} h(x,t) \, \mathrm dt = s^{\gamma+\nu-1}
				\left( 1+s^{-\nu} \right)^\delta e^{-x s^{\gamma+\nu}\left( 1+s^{-\nu} \right)^\delta},
				\qquad x \ge 0, \: s>0, \\
				& \label{05lap2} \tilde{h}'(z,t) = \int_0^\infty e^{-zx} h(x,t) \, \mathrm dx
				=\sum_{r=0}^\infty (-z)^r t^{r(\gamma+\nu)} E_{\nu,r(\gamma+\nu)+1}^{r\delta} (-t^\nu),
				\qquad t \ge 0, \: z > 0,
			\end{align}
			where $|z/(s^{\gamma+\nu}(1+ s^{-\nu})^\delta)|<1$.
			
			\begin{proof}
				Formula \eqref{05lap} is straightforward as it follows from the expression of the
				Laplace transform of an exponential. For what concerns formula \eqref{05lap2},
				from \eqref{hlaplap}, we first write
				\begin{align}
					\tilde{\tilde{h}}(z,s) = s^{-1} \left( 1+z s^{-(\gamma+\nu)} \left( 1+s^{-\nu} \right)^{-\delta}
					\right)^{-1}
					= \sum_{r=0}^\infty (-z)^r s^{-(\gamma+\nu)r-1} \left( 1+s^{-\nu} \right)^{-r\delta},
				\end{align}
				which holds if $|z/(s^{\gamma+\nu}(1+ s^{-\nu})^\delta)|<1$.			
				We then recall the following formula for the Laplace transform of a generalized
				Mittag--Leffler function \citep[formula (2.3.24)]{mathai}
				\begin{align}
					\label{20de}
					\int_0^\infty t^{\eta-1} e^{-pt} E_{\alpha, \eta}^\xi (\zeta t^\alpha)\, \mathrm dt =
					p^{-\eta} \left( 1-\zeta p^{-\eta} \right)^{-\xi}, \qquad \alpha>0, \: \eta>0,
					\: p > |\zeta|^{1/\alpha},
				\end{align}
				by means of which result \eqref{05lap2} is easily found.
			\end{proof}
		\end{proposition}		
		
		\begin{remark}
			Note how the above $t$-Laplace transform \eqref{05lap}, for $\delta=0$, reduces to the well-known
			$t$-Laplace transform of the inverse $(\gamma+\nu)$-stable subordinator, namely
			$s^{\gamma+\nu-1} \exp (-xs^{\gamma+\nu})$ \citep[formula (2.14)]{mirko}.
			Furthermore, the $x$-Laplace transform, again in the case $\delta=0$, can be written as
			\begin{align}
				\tilde{h}'(z,t) = \sum_{r=0}^\infty \frac{(-z t^{\gamma+\nu})^r}{\Gamma(r(\gamma+\nu)+1)} =
				E_{\gamma+\nu,1}(-z t^{\gamma+\nu}),
			\end{align}
			which, as expected, concides with the classical result \citep[formula (2.13)]{mirko}.
		\end{remark}

		The remaining pages of this section are devoted to explain in which sense the
		integral operator we are analyzing is connected to some stochastic processes. For the sake of
		clarity we explain our results first in the specific case $\delta \in (0,1]$ leaving the
		presentation of the more general case $\delta > 0$ at the end of this section.
	
		Consider a filtered probability space $(\Omega, \mathcal{F}, \mathfrak{G}, \mathbb{P})$,
		where $\mathfrak{G}=(\mathcal{G}_t)_{t \ge 0}$ is the associated filtration, and		
		the process
		\begin{align}
			U_t^{(\alpha_1, \alpha_2)}
			= {}_1 V_t^{\alpha_1} + {}_2 V_t^{\alpha_2}, \qquad t \ge 0, \: \alpha_1,\alpha_2 \in (0,1),
		\end{align}
		adapted to $\mathfrak{G}$,
		where ${}_jV_t^{\alpha_j}$, $t \ge0$, $j=1,2$, are independent stable subordinators of order $\alpha_j$. Let us
		further consider the stable
		subordinator $V^\delta_t$, $t \ge 0$, $\delta \in (0,1]$ also adapted to $\mathfrak{G}$ and independent of
		$U_t^{(\alpha_1, \alpha_2)}$, $t \ge 0$.
		We focus now on the subordinated process
		\begin{equation}
			\label{Vproc}
			U_{V^\delta_t}^{(\alpha_1, \alpha_2)} = {}_1 V_{V^\delta_t}^{\alpha_1}
			+ {}_2V_{V^\delta_t}^{\alpha_2},
			\qquad t \ge 0,
		\end{equation}
		clearly adapted to the time-changed filtration $\bigl( \mathcal{G}_{V^\delta_t}\bigr)_{t \ge 0}$
		and, in particular we have that its Laplace transform is
		\begin{align}
			\label{LapVproc}
			\mathbb{E} \exp \left( -z U_{V^\delta_t}^{(\alpha_1, \alpha_2)} \right)
			= \mathbb{E}
			\exp \left( - z^{\alpha_1} V^\delta_t - z^{\alpha_2} V^\delta_t \right)
			= \exp \left( - t (z^{\alpha_1} + z^{\alpha_2})^\delta \right).
		\end{align}
		Notice that for \eqref{Vproc} we obtain
		\begin{equation}
			\label{VVproc}
			U_{V^\delta_t}^{(\alpha_1, \alpha_2)}  \stackrel{\mathrm d}{=} {}_1V_t^{\delta \alpha_1}
			+ {}_2V_t^{\delta \alpha_2}, \qquad t \ge 0,
		\end{equation}
		where ${}_j V_t^{\delta \alpha_j}$, $t \ge 0$, $j=1,2$ are now $\mathfrak{G}$-adapted
		dependent stable subordinators. The dependence is due to the time-change and vanishes
		in the degenerate case $\delta=1$.
		For $\delta \in (0,1)$ the processes in \eqref{VVproc} possess dependent increments but non-decreasing paths,
		that is the increments are non-negative. Thus we argue that
		\begin{align}
			\label{Eproc}
			C_t^{(\delta\alpha_1, \delta\alpha_2)} = \inf\left\lbrace s\geq 0 \colon {}_1 V_s^{\delta \alpha_1}
			+ {}_2V_s^{\delta \alpha_2} \notin (0, t) \right\rbrace
            \stackrel{\mathrm d}{=} \inf \Bigl\{ s\geq 0
			\colon	U_{V^\delta_s}^{(\alpha_1, \alpha_2)} \notin (0, t) \Bigr\}, \qquad t \ge 0,
		\end{align}
		is the first exit time of the process \eqref{VVproc} from the interval $(0,t)$
		whose distribution coincides with that
		of the first exit time of \eqref{Vproc} from the same interval $(0,t)$. We refer to \eqref{Eproc}
		also as the inverse to \eqref{Vproc} as
		\begin{equation}
			\mathbb{P} \left\lbrace C_t^{(\delta \alpha_1, \delta \alpha_2)} > x \right\rbrace
			= \mathbb{P} \left\lbrace U_{V^\delta_x}^{(\alpha_1, \alpha_2)} < t \right\rbrace.
		\end{equation}
		Note that if we let the inverse process $C_t^{(\delta\alpha_1, \delta\alpha_2)}$, $t \ge 0$,
		be adapted to filtration $\mathfrak{F} = (\mathcal{F}_t)$, we have that $\mathfrak{G} = (\mathcal{G}_t)_{t\ge 0}
		= \bigl(\mathcal{F}_{C_t^{(\delta\alpha_1, \delta\alpha_2)}}\bigr)_{t \ge 0}$, that is
		$C_t^{(\delta\alpha_1, \delta\alpha_2)}$ is a time change on the filtered probability
		space $(\Omega, \mathcal{F}, \mathfrak{F}, \mathbb{P})$.
		Relation \eqref{Eproc} allows us to derive the $t$-Laplace transform of the density law of \eqref{Eproc}: 
		\begin{align}
			- \frac{\mathrm d}{\mathrm d x} & \int_0^\infty e^{-s t}
			\mathbb{P} \left\lbrace C_t^{(\delta \alpha_1, \delta \alpha_2)} > x
			\right\rbrace \mathrm dt = -\frac{\mathrm d}{\mathrm d x} \int_0^\infty e^{-s t} \mathbb{P} \left\lbrace
			U_{V^\delta_x}^{(\alpha_1, \alpha_2)} < t \right\rbrace  \mathrm dt \\
			= & - \frac{1}{s}\frac{\mathrm d}{\mathrm d x} \int_0^\infty e^{-s t}
			\mathbb{P} \left\lbrace U_{V^\delta_x}^{(\alpha_1,
			\alpha_2)} \in \mathrm dt \right\rbrace
			= - \frac{1}{s}\frac{\mathrm d}{\mathrm d x} \mathbb{E}
			\exp \left( -s U_{V^\delta_x}^{(\alpha_1, \alpha_2)} \right). \notag
		\end{align}
		Now, by using the Laplace transform \eqref{LapVproc} we arrive at
		\begin{align}
			\label{tLapEproc}
			\int_0^\infty e^{-s t}  \left( \mathbb{P} \left\lbrace C_t^{(\delta \alpha_1, \delta \alpha_2)}
			\in \mathrm dx 
			\right\rbrace / \mathrm dx\right) \mathrm dt = & - \frac{1}{s} \frac{\mathrm d}{\mathrm d x} \exp
			\left( - x (s^{\alpha_1} + s^{\alpha_2})^\delta \right) \\
			= & \frac{1}{s} (s^{\alpha_1} + s^{\alpha_2})^\delta \exp
			\left( - x (s^{\alpha_1} + s^{\alpha_2})^\delta \right). \notag 
		\end{align}
	
		First we present the following result.

		\begin{proposition}
			\label{Prop-L-K}
			We have that
			\begin{equation}
				C_t^{(\delta \alpha_1, \delta \alpha_2)} \overset{\mathrm d}{=}
				L^\delta_{C_t^{(\alpha_1, \alpha_2)}}, \qquad t \ge 0, \: \delta \in (0,1],
			\end{equation}
			where $L^\delta_t = \inf\{x \geq 0\, \colon \, V^\delta_x \notin (0,t)\}$,
			$t \ge 0$ is the inverse to the stable subordinator $V^\delta_t$, $t \ge 0$ in the sense that
			\begin{align}
				\mathbb{P} \{ L^\delta_t < x \} = \mathbb{P} \{V^\delta_x > t \}
			\end{align}
			and $C_t^{(\alpha_1, \alpha_2)} = \inf\left\lbrace s\geq 0
			\colon {}_1 V_s^{ \alpha_1} + {}_2 V_s^{ \alpha_2} \notin (0, t) \right\rbrace$, $t \ge 0$,
			is the inverse to $U_{t}^{(\alpha_1, \alpha_2)}$, $t \ge 0$ in the sense that
			\begin{equation}
				\mathbb{P} \left\lbrace C_t^{(\alpha_1, \alpha_2)} > x \right\rbrace
				= \mathbb{P} \left\lbrace {}_1V_x^{\alpha_1} + {}_2V_x^{\alpha_2} < t \right\rbrace.
			\end{equation}
	
			\begin{proof}
				It suffices to consider formula \eqref{tLapEproc} for $\delta = 1$ and the integral
				\begin{align}
					\int_0^\infty & e^{-st} \Big( \mathbb{E}
					\exp\big( - z L^\delta_{C_t^{(\alpha_1, \alpha_2)}} \big)
					\Big) \mathrm dt = \int_0^\infty \mathbb{E} \left[e^{-z L^\delta_x}\right] \int_0^\infty e^{-s t}
					\mathbb{P} \left\lbrace C_t^{ (\alpha_1, \alpha_2)}\in \mathrm dx
					\right\rbrace \mathrm dt \\
					= & \int_0^\infty \mathbb{E} \left[e^{-z L^\delta_x}\right] \frac{1}{s} (s^{\alpha_1} + s^{\alpha_2})
					\exp\left( - x (s^{\alpha_1} + s^{\alpha_2}) \right) \mathrm dx \notag \\
					= & \frac{1}{s} (s^{\alpha_1} + s^{\alpha_2}) \int_0^\infty E_\delta(-z x^\delta)
					\, \exp\left( - x (s^{\alpha_1} + s^{\alpha_2}) \right) \mathrm dx
					= \frac{1}{s}  \frac{(s^{\alpha_1} + s^{\alpha_2})^\delta}{z
					+ (s^{\alpha_1} + s^{\alpha_2})^\delta} \notag
				\end{align}
				which coincides with the $x$-Laplace transform of \eqref{tLapEproc}. 
			\end{proof}
		\end{proposition}
		
		In the following, when we refer to stochastic solution of a pde, we mean the stochastic
		process whose density function is the fundamental solution to such pde.

		\begin{remark}
			We observe that the process $C^{(\alpha_1, \alpha_2)}_t$, $t \ge 0$ has been investigated in
			\citet{DTO2012} and is, for $\alpha_1 = \gamma+\nu \in (0,1)$ and $\alpha_2 = \gamma \in (0,1)$
			the stochastic solution to the fractional telegraph equation 
			\begin{equation}
				\label{cauchy-prob-K}
				\frac{\partial^{\gamma+\nu}}{\partial t^{\gamma+\nu}} u(x,t)
				+ \frac{\partial^\gamma}{\partial t^\nu} u(x,t)
				= - \frac{\partial}{\partial x} u(x,t), \qquad x \ge 0, \: t \ge 0,
			\end{equation}
			subject to the initial and the boundary conditions
			\begin{align}
				u(x, 0) = \delta(x), \qquad
				u(0, t) = \frac{t^{-\gamma -\nu }}{\Gamma(1-\gamma-\nu)} + \frac{t^{-\nu}}{\Gamma(1-\nu)}.
			\end{align}
		\end{remark}

		Now we are ready to prove the following result which shows the relation to the Cauchy problem \eqref{04sub}.

		\begin{theorem}
			\label{gardel}
			The stochastic solution to
			\eqref{04sub}, $\delta \in (0,1]$
			is given by the hitting time \eqref{Eproc} of the subordinated process \eqref{Vproc} with 
			$\alpha_1 = (\gamma + \nu)/\delta \in (0,1]$, $\alpha_2	= (\gamma + \nu)/\delta - \nu \in (0,1]$.
			Furthermore, the process \eqref{Eproc} becomes
			\begin{align}
				Z_t^\delta := C_t^{(\gamma + \nu, \gamma + \nu - \delta \nu)}
				= \inf\left\lbrace s\geq 0 \colon
				{}_1 V_s^{\gamma + \nu} + {}_2 V_s^{\gamma + \nu - \delta \nu}
				\notin (0, t) \right\rbrace, \qquad t \ge 0,
			\end{align}
			where $ {}_1 V_t^{\gamma + \nu}$ and  ${}_2 V_t^{\gamma + \nu - \delta \nu}$
			are dependent stable subordinators.
			\begin{proof}
				From \eqref{tLapEproc}
				we obtain the Laplace--Laplace transform of the density law of the process
				\eqref{Eproc} as follows
				\begin{align}
					\int_0^\infty & e^{-st} \, \mathbb{E}\exp
					\left( -z C_t^{(\delta \alpha_1, \delta \alpha_2)} \right) \mathrm dt
					= \int_0^\infty \int_0^\infty e^{-s t- z x}
					\mathbb{P} \left\lbrace C_t^{(\delta \alpha_1,
					\delta \alpha_2)} \in \mathrm dx \right\rbrace \mathrm dt \\
					= & \int_0^\infty e^{-z x} \frac{1}{s} (s^{\alpha_1}
					+ s^{\alpha_2})^\delta \exp \left( - x (s^{\alpha_1}
					+ s^{\alpha_2})^\delta \right) \mathrm dx
					= \frac{1}{s}  \frac{(s^{\alpha_1} + s^{\alpha_2})^\delta}{z
					+ (s^{\alpha_1} + s^{\alpha_2})^\delta}. \notag
				\end{align}
				For $\alpha_1 = (\gamma + \nu)/\delta$, $\alpha_2 = (\gamma + \nu)/\delta-\nu$,
				we have that
				\begin{align}
					\int_0^\infty e^{-s t} \, \mathbb{E} \exp\left( - z
					Z_t^\delta
					\right) \mathrm dt
					= & \frac{s^{\gamma + \nu -1} (1+ s^{-\nu})^\delta}{z + s^{\gamma + \nu} (1+ s^{-\nu})^\delta},
				\end{align}
				which coincides with \eqref{hlaplap}.
			\end{proof}
		\end{theorem}
		
		We now move to analyzing the more general case $\delta > 0$.
		In light of Theorem \ref{gardel}, the results stated in the following theorem will appear natural.
		What changes is basically that now we are dealing with a linear combination of subordinated
		stable subordinators whose hitting time will be the stochastic solution to
		\eqref{04sub}.
		
		\begin{theorem}
			\label{55teo}
			Given the filtered probability space $(\Omega, \mathcal{F}, \mathfrak{F}, \mathbb{P})$,
			the stochastic solution to
			\eqref{04sub}, $\delta > 0$, is given by the $\mathfrak{F}$-hitting time $Z_t^\delta$, $t \ge 0$,
			of the $\mathfrak{G}$-adapted process
			\begin{align}
				\label{ddir}
				M^\delta_t & = \sum_{r=0}^n \binom{n}{r}^{1/[(\gamma+\nu)n/\delta-r\nu]}
				{}_r V_t^{\gamma+\nu -r \nu \delta/n} =
				\sum_{r=0}^n {}_r V_{\binom{n}{r}t}^{\gamma+\nu -r \nu \delta/n}, \qquad t \ge 0,
			\end{align}
			where $\mathfrak{F} = (\mathcal{F}_t)_{t \ge 0}
			= \bigl( \mathcal{G}_{Z_t^\delta} \bigr)_{t \ge 0}$,
			is the associated filtration (with $\mathfrak{G}=(\mathcal{G}_t)_{t \ge 0}$),
			${}_r V_t^{\gamma+\nu -r \nu \delta/n}$,
			$r=1,\dots n$, are dependent stable subordinators,
			$n= \lceil \delta \rceil$ is the ceiling of $\delta$ and $\nu \delta < \gamma+\nu < 1$.
			
			\begin{proof}
				We start, similarly to proof of Theorem \ref{gardel}, by considering the process
				\begin{align}
					U_t & = \sum_{r=0}^n \binom{n}{r}^{1/[(\gamma+\nu)n/\delta-r\nu]}
					{}_r V_t^{(\gamma+\nu)n/\delta-r\nu},
				\end{align}
				with ${}_r V_t^{(\gamma+\nu)n/\delta-r\nu}$, $r=1,\dots n$, independent stable subordinators
				and $\nu \delta < \gamma+\nu < \delta/n$.
				This process must be subordinated to a further stable subordinator $V_t^{\delta/n}$ independent of
				${}_r V_t^{(\gamma+\nu)n/\delta-r\nu}$, $r=1,\dots n$. We thus obtain
				\begin{align}
					U_{V_t^{\delta/n}} = \sum_{r=0}^n \binom{n}{r}^{1/[(\gamma+\nu)n/\delta-r\nu]}
					{}_r V_{V_t^{\delta/n}}^{(\gamma+\nu)n/\delta-r\nu}.
				\end{align}
				Note that, as $U_t$ is a linear combination of independent stable subordinators, the process
				$U_{V_t^{\delta/n}} \overset{\mathrm d}{=} M^\delta_t$.
				Due to the nature of the time-change considered it is clear that $M_t^\delta$
				is a linear combination of dependent subordinators. This is easily explained by noticing that the shared time-change
				turns the independency $({}_0 V_t^{(\gamma+\nu)n/\delta}, {}_1 V_t^{(\gamma+\nu)n/\delta-1\nu},{}_2 V_t^{(\gamma+\nu)2/\delta-\nu}
				\dots, {}_n V_t^{(\gamma+\nu)n/\delta-n\nu})$ into the dependent collection
				of subordinators $({}_0 V_t^{\gamma+\nu}, {}_1 V_t^{\gamma+\nu -\nu \delta/n},
				{}_2 V_t^{\gamma+\nu -2 \nu \delta/n}, \dots, {}_n V_t^{\gamma+\nu -n \nu \delta/n})$.
				The Laplace transform of the density law of $U_{V_t^{\delta/n}}$
				reads
				\begin{align}
					\mathbb{E} \exp \left\{-z U_{V_t^{\delta/n}}\right\}
					& = \mathbb{E} \exp \Biggl\{ - V_t^{\delta/n} \sum_{r=0}^n
					\binom{n}{r} z^{(\gamma+\nu)n/\delta-r\nu} \Biggr\}
					= \exp \Biggl\{ -t \left[ \sum_{r=0}^n \binom{n}{r} z^{(\gamma+\nu)n/\delta-r\nu}
					\right]^{\delta/n} \Biggr\}  \\
					& = \exp \left\{ -t \left[ z^{(\gamma+\nu)n/\delta} (1+z^{-\nu})^n \right]^{\delta/n} \right\}
					\notag
					= \exp \left\{ -t \, z^{\gamma+\nu} (1+z^{-\nu})^\delta \right\}. \notag
				\end{align}
				
				Let us define now the right-inverse process to $M^\delta_t$ as
				\begin{align}
					Z_t^\delta := \inf \left\{ s \ge 0 \colon M^\delta_t \notin (0,t) \right\}
					\overset{\mathrm d}{=} \inf \left\{ s \ge 0 \colon U_{V_t^{\delta/n}} \notin (0,t) \right\},
					\qquad t \ge 0, \: \delta >0.
				\end{align}
				In particular note that $\mathbb{P} \{ Z_t^\delta >x \}
				= \mathbb{P} \{ U_{V_x^{\delta/n}} < t \}$.
				
				The time-Laplace transform related to the inverse process $Z_t^\delta$, $t \ge 0$,
				can be determined as in the following.
				\begin{align}
					-\frac{\mathrm d}{\mathrm dx} & \int_0^\infty e^{-s t} \mathbb{P}
					\{ Z_t^\delta > x \} \mathrm dt
					= -\frac{\mathrm d}{\mathrm dx} \int_0^\infty e^{-st} \mathbb{P}
					\{ U_{V_x^{\delta/n}} < t \}
					\mathrm dt \\
					& = - s^{-1} \frac{\mathrm d}{\mathrm dx} \int_0^\infty e^{-st}
					\mathbb{P} \{ U_{V_x^{\delta/n}} \in \mathrm dt \}
					= - s^{-1} \frac{\mathrm d}{\mathrm dx} \mathbb{E} e^{-s U_{V_x^{\delta/n}}}. \notag
				\end{align}
				Therefore
				\begin{align}
					\label{tetration}
					\int_0^\infty e^{-st} \left[ \mathbb{P} \{ Z_t^\delta
					\in \mathrm dx \} / \mathrm dx \right]
					\mathrm dt & = - s^{-1} \frac{\mathrm d}{\mathrm dx}
					\exp \left\{ -x \, s^{\gamma+\nu} (1+s^{-\nu})^\delta \right\} \\
					& = s^{\gamma+\nu-1} (1+s^{-\nu})^\delta \exp \left\{ -x \, s^{\gamma+\nu}
					(1+s^{-\nu})^\delta \right\}. \notag
				\end{align}
				
				Finally we calculate the complete Laplace--Laplace transform.
				\begin{align}
					\int_0^\infty & e^{-st} \mathbb{E} e^{-z Z_t^\delta} \mathrm dt =
					\int_0^\infty \int_0^\infty e^{-st-zx} \mathbb{P}
					\{ Z_t^\delta \in \mathrm dx \} \mathrm dt \\
					& = \int_0^\infty e^{-zx} 
					s^{\gamma+\nu-1} (1+s^{-\nu})^\delta \exp \left\{ -x \, s^{\gamma+\nu}
					(1+s^{-\nu})^\delta \right\} \mathrm dx
					= \frac{s^{\gamma+\nu-1} (1+s^{-\nu})^\delta}{z + s^{\gamma+\nu} (1 + s^{-\nu})^\delta}. \notag
				\end{align}
				As the latter expression coincides with \eqref{hlaplap} the proof of the theorem is complete.
			\end{proof}
		\end{theorem}
		
		\begin{remark}
			If we consider $n=\lceil \delta \rceil$ and  $m > n$ we have that
			the corresponding process becomes
			\begin{align}
				M^\delta_{V_t^{n/m}} = \sum_{r=0}^n \binom{n}{r}^{1/[(\gamma+\nu)n/\delta-r\nu]}
				{}_r V_t^{(\gamma+\nu)n/m - r \nu \delta/m}, \qquad t \ge 0, \: \forall m > n,
			\end{align}
			where $\nu \delta < \gamma+\nu < 1$ and $V_t^{n/m}$ is an independent $n/m$-stable
			subordinator.
			It is worthy noticing that in this case the hitting time of the above process is not a stochastic solution
			to equation \eqref{04sub}.
		\end{remark}
		
		Aside the results contained in Theorem \ref{55teo} we are able also to prove a subordination relation
		for the stochastic solution $Z_t^\delta$, $t \ge 0$, in the following proposition.
		
		\begin{proposition}
			\label{mariolino}
			We have that
			\begin{align}
				Z^\delta_t \overset{\mathrm d}{=} L^{\delta/n}_{F_t}, \qquad t \ge 0,
			\end{align}
			where $n = \lceil \delta \rceil$ is the ceiling of $\delta > 0$ and
			where $L^{\delta/n}_t = \inf\{x \geq 0\, \colon \, V^{\delta/n}_x \notin (0,t)\}$,
			$t \ge 0$, is the right-inverse process to the stable subordinator $V^{\delta/n}_t$, $t \ge 0$,
			in the sense that
			$\mathbb{P} \{ L^{\delta/n}_t < x \} = \mathbb{P} \{V^{\delta/n}_x > t \}$
			and the hitting time $F_t = \inf\left\lbrace s\geq 0
			\colon U_s \notin (0, t) \right\rbrace$, $t \ge 0$
			is the inverse to $ U_t$, $t \ge 0$ in the sense that
			$\mathbb{P} \left\lbrace F_t > x \right\rbrace
				= \mathbb{P} \left\lbrace U_x < t \right\rbrace$.
			
			\begin{proof}
				In order to prove the subordination relation
				it is sufficient to consider formula \eqref{tetration} for $\delta = n$
				and the following calculations.
				\begin{align}
					& \int_0^\infty e^{-st} \Big( \mathbb{E} \exp\big( - z L^{\delta/n}_{F_t } \big)
					\Big) \mathrm dt
                    = \int_0^\infty \mathbb{E}
					\left[e^{-z L^{\delta/n}_x}\right] \int_0^\infty e^{-s t}
					\mathbb{P} \left\lbrace F_t \in \mathrm dx \right\rbrace  \, \mathrm dt \\
					& = \int_0^\infty \mathbb{E} \left[e^{-z L^{\delta/n}_x}\right]
					s^{(\gamma+\nu)n/\delta-1} (1+s^{-\nu})^n
					\exp\left( - x \, s^{(\gamma + \nu)n/\delta} (1+s^{-\nu})^n \right) \mathrm dx \notag \\
					& = s^{(\gamma+\nu)n/\delta-1} (1+s^{-\nu})^n \int_0^\infty E_{\delta/n, 1} (-z x^{\delta/n})
					\exp\left( - x \, s^{(\gamma + \nu)n/\delta} (1+s^{-\nu})^n \right) \mathrm dx \notag \\
					& = s^{(\gamma+\nu)n/\delta-1} (1+s^{-\nu})^n \frac{\left[ s^{(\gamma+\nu)n/\delta}
					(1+s^{-\nu})^n \right]^{\delta/n-1}}{
					z + \left[ s^{(\gamma+\nu)n/\delta}(1+s^{-\nu})^n \right]^{\delta/n}}
					= \frac{s^{\gamma+\nu-1}(1+s^{-\nu})^\delta}{z + s^{\gamma+\nu}(1+s^{-\nu})^\delta}. \notag
				\end{align}
				The last expression exactly coincides with the $x$-Laplace transform of \eqref{tetration}. Notice that
				$E_{\delta/n,1}(x) = E_{\delta/n,1}^1(x)$ is the classical two-parameter Mittag--Leffler function.
			\end{proof}
		\end{proposition}

		\begin{remark}
			In the specific case of $\delta=n \in \mathbb{N}\cup \{ 0 \}$, equation \eqref{04sub}
			takes a peculiar form. We have
			\begin{align}
				\left( \mathbb{D}_{\nu,\gamma+\nu,-1;0+}^n h(x,\cdot) \right) (t)
				= - \frac{\partial}{\partial x} h(x,t) 
				\Leftrightarrow \: &
				\frac{\partial^{\gamma+\nu+\theta}}{\partial t^{\gamma+\nu+\theta}}
				\int_0^t (t-y)^{\theta-1} E_{\nu,\theta}^{-n} \left[ - (t-y)^\nu \right]
				h(x,y) \, \mathrm dy \\
				& \qquad = - \frac{\partial}{\partial x} h(x,t) + \delta(x) \, t^{-(\gamma+\nu)}
				E_{\nu,1-(\gamma+\nu)}^{-n} (- t^\nu). \notag
			\end{align}
			By recalling again that $(-n)_r = (-1)^r (n-r+1)_r = (-1)^r n!/(n-r)!$,
			we obtain
			\begin{align}
				\label{30tt}
				& \frac{\partial^{\gamma+\nu+\theta}}{\partial t^{\gamma+\nu+\theta}}
				\int_0^t (t-y)^{\theta-1} \sum_{r=0}^n \binom{n}{r} \frac{(t-y)^{\nu r}}{\Gamma(\nu r+\theta)}
				h(x,y) \, \mathrm dy
				= - \frac{\partial}{\partial x} h(x,t) + \delta(x)
				\sum_{r=0}^n \binom{n}{r} \frac{t^{-\left(\gamma-\nu\left(r-1\right)\right)}}{\Gamma \left(1-
				\left(\gamma-\nu\left(r-1\right)\right)\right)} \\
				& \Leftrightarrow \quad
				\sum_{r=0}^n \binom{n}{r} \frac{\partial^{\gamma-\nu(r-1)}}{\partial
				t^{\gamma-\nu(r-1)}} h(x,t)
				= - \frac{\partial}{\partial x} h(x,t) + \delta(x)
				\sum_{r=0}^n \binom{n}{r} \frac{t^{-\left(\gamma-\nu\left(r-1\right)\right)}}{\Gamma \left(1-
				\left(\gamma-\nu\left(r-1\right)\right)\right)} \notag \\
				& \Leftrightarrow \quad
				\sum_{r=0}^n \binom{n}{r} \frac{\mathfrak{d}^{\gamma-\nu(r-1)}}{\mathfrak{d}
				t^{\gamma-\nu(r-1)}} h(x,t)
				= - \frac{\partial}{\partial x} h(x,t) \notag
			\end{align}
			with $h(x,0)=\delta(x)$, $0 < \gamma-\nu(r-1) < 1$ and thus $n\nu < \gamma+\nu < 1$.
			Note furthermore that equations \eqref{30tt} and \eqref{22tt} are consistent
			with Theorem 3.1 of \citet{umarov}.
		\end{remark}
		
		Before moving to the next section where the introduced process is related to the
		the stochastic solution of different abstract Cauchy problems,
		we first underline in the following remark that the inverse process $Z_t^\delta$, $t \ge 0$,
		$\delta > 0$, is well-behaved and can be used as a time-change.
		
		\begin{remark}
			The inverse process $Z_t^\delta$, $t \ge 0$,
			$\delta > 0$, is a continuous time-change on the probability space
			$(\Omega, \mathcal{F}, \mathfrak{F}, \mathbb{P})$.
			This simply ensues from the construction of the process \eqref{ddir}
			as a linear combination  with non-negative coefficients
			of dependent stable subordinators which are clearly increasing processes
			(right-continuous and taking values in $[0,\infty]$) and $M_t^\delta$, $t \ge 0$,
			is adapted to $\mathfrak{G}$.
		\end{remark}
		
		We conclude this section by studying the case $\delta<0$. When $\delta$ is strictly
		negative calculations become more complicated. We present the following result for $-1/2<\delta<0$,
		$\nu=-\beta$.

		\begin{theorem}

			Let $0 < \gamma - \beta < 1$, (with $\gamma, \beta \in (0,1)$)
			and $\epsilon = -\delta$ with $ -1/2 < \delta < 0$. In this case the solution to \eqref{04sub} is
			\begin{equation}
				h(x,t) = \frac{1}{\sqrt{\pi}}  \int_0^\infty \frac{\mathrm dz}{\sqrt{z}} \phi(x,z)
				\int_0^\infty \mathrm du \, \varphi(u,z,t), \quad x,t >0,
			\end{equation}
			where
			\begin{equation}
				\varphi(u,z,t) =  \int_0^t \int_0^\infty e^{-y} \, v_{2\epsilon}(y, u^2/4z)
				\, v_{\beta}(t-k, y^{\beta})\, l_{\gamma - \beta}(u,k)\, \mathrm dy\, \mathrm dk.
			\end{equation}
			$v_a(x,t) = \mathbb{P} \{ V_t^a \in \mathrm dx \}/\mathrm dx$ is the density law of
			the $a$-stable subordinator, $l_a(x,t) = \mathbb{P} \{ L_t^a \in \mathrm dx \}/\mathrm dx$
			is the density law of the inverse process to an $a$-stable subordinator,
			and $\phi$ is a function such that
			$\tilde{\phi}(\mu,z) = \mu e^{-z \mu^2}$.

			\begin{proof}
				First we show that
				\begin{equation}
					\label{cappuccino}
					\tilde{\varphi}(u,z,s) = e^{- \frac{u^2}{4z} (1+s^\beta)^{2\epsilon}}
					s^{\gamma-\beta-1} e^{- u s^{\gamma -\beta}}.
				\end{equation}
				Indeed, the first exponential term in the right-hand side of \eqref{cappuccino} can be written as
				\begin{align}
					\label{bloc}
					e^{- \frac{u^2}{4z} (1+s^\beta)^{2\epsilon}}
					= & \mathbb{E} \, e^{-(1+s^\beta) V^{2\epsilon}_{u^2/4z}}
					= \mathbb{E} \, \mathbb{E} \, e^{- (1 + s V^\beta_1)
					V^{2\epsilon}_{u^2/4z}} \\
					= & \int_0^\infty e^{-y} \mathbb{E} \, e^{-ys V^\beta_1} v_{2\epsilon}(y, u^2/4z)
					\, \mathrm dy
					= \int_0^\infty e^{-y -(ys)^\beta} v_{2\epsilon}(y, u^2/4z) \, \mathrm dy, \notag
				\end{align}
				where
				$e^{-(ys)^\beta} = \int_0^\infty e^{-s x} v_{\beta}(x, y^{\beta}) \, \mathrm dx$.
				Also, the remaining terms of equation \eqref{cappuccino} are in fact the Laplace transform
				of the density law of an inverse $(\gamma-\beta)$-stable subordinator, that is
				\begin{align}
					\label{babaru}
					s^{\gamma-\beta-1} e^{- u s^{\gamma -\beta}} = \int_0^\infty
					e^{-s x} l_{\gamma - \beta}(u, x)\, \mathrm dx.
				\end{align}
				Now, note that the following product can be represented
				as the Laplace transform of a convolution, as it is shown in the following formula. 
				\begin{align}
					\label{mareblu}
					e^{-(ys)^\beta} s^{\gamma-\beta-1} e^{- u s^{\gamma -\beta}}
					= \int_{0}^\infty e^{-s t} \int_0^t v_{\beta}(t-w, y^{\beta}) l_{\gamma - \beta}(u, w) \, \mathrm dw
					\, \mathrm dt.
				\end{align}
				Therefore we have
				\begin{equation}
					\int_{0}^\infty e^{-s t}\varphi(u,z,t)\, \mathrm dt
					= \int_0^\infty e^{-y} v_{2 \epsilon} (y,u^2/4z)
					e^{-(ys)^\beta} s^{\gamma-\beta-1} e^{- u s^{\gamma -\beta}} \mathrm dy.
				\end{equation}
				Last step has been obtained by using \eqref{mareblu}. Considering \eqref{bloc}
				we immediately obtain result \eqref{cappuccino}.

				The $t$-Laplace transform of $h$ becomes therefore
				\begin{equation}
					\tilde{h}(x,s) = \frac{1}{\sqrt{\pi}} \int_0^\infty \frac{\mathrm dz}{\sqrt{z}} \phi(x,z)
					\int_0^\infty e^{- \frac{u^2 (1+s^\beta)^{2\epsilon}}{4z}}
					s^{\gamma-\beta-1} e^{- u s^{\gamma -\beta}} \, \mathrm du
				\end{equation}
				and the double Laplace transform is given by
				\begin{equation}
					\tilde{\tilde{h}}(\mu,s) = \frac{1}{\sqrt{\pi}} \int_0^\infty \frac{\mathrm dz}{\sqrt{z}}
					\left( \mu e^{-z \mu^2} \right) \int_0^\infty e^{- \frac{u^2 (1+s^\beta)^{2\epsilon}}{4z}}
					s^{\gamma-\beta-1} e^{- u s^{\gamma -\beta}} \, \mathrm du.
				\end{equation}
				By considering that
				$K_\alpha(2\sqrt{ab}) = \frac{1}{2} \left( \frac{b}{a} \right)^{\frac{\alpha}{2}}
				\int_0^\infty y^{\alpha -1} e^{-yb - \frac{a}{y}} \, \mathrm dy$
				is the modified Bessel function of the second kind, we have that
				\begin{align}
					\int_0^\infty \frac{\mathrm dz}{\sqrt{z}} \left( \mu e^{-z \mu^2} \right)
					e^{- \frac{u^2 (1+s^\beta)^{2\epsilon}}{4z}}
					= 2 \mu \left(\frac{u^2 (1+s^\beta)^{2\epsilon}}{4 \mu^2} \right)^\frac{1}{4} K_\frac{1}{2}
					\left( \sqrt{\mu^2 u^2 (1+s^\beta)^{2\epsilon}} \right).
				\end{align}
				Since $K_{1/2}(z) = \sqrt{\frac{\pi}{z}} e^{-z}$, we get
				\begin{align}
					\int_0^\infty \frac{\mathrm dz}{\sqrt{z}} \left( \mu e^{-z \mu^2} \right)
					e^{- \frac{u^2 (1+s^\beta)^{2\epsilon}}{4z}}
					= 2 \mu \left(\frac{u^2 (1+s^\beta)^{2\epsilon}}{4 \mu^2} \right)^\frac{1}{4}
					\sqrt{\frac{\pi}{2}} \left( \mu^2 u^2 (1+s^\beta)^{2\epsilon} \right)^{-\frac{1}{4}}
					e^{- \mu u (1+s^\beta)^{\epsilon}}
					= \sqrt{\pi} e^{- \mu u (1+s^\beta)^{\epsilon}}.
				\end{align}
				Hence, 
				\begin{align}
					\tilde{\tilde{h}}(\mu,s) = s^{\gamma-\beta-1} \int_0^\infty \mathrm du \, e^{-u s^{\gamma -\beta}
					- \mu u (1+s^\beta)^{\epsilon}}
					= \frac{s^{\gamma-\beta-1}}{s^{\gamma -\beta} + \mu  (1+s^\beta)^{\epsilon}}
					= \frac{s^{\gamma-\beta-1}(1+s^\beta)^{-\epsilon}}{s^{\gamma -\beta} (1+s^\beta)^{-\epsilon} + \mu }
				\end{align}
				which coincides with \eqref{hlaplap} for $\beta=-\nu$ and $\epsilon=-\delta$.
			\end{proof}
		\end{theorem}

	\section{Time changed L\'evy processes}
	
        \label{lele}
		Let $A^x_t$, $t \ge 0$, be an $\mathbb{R}^d$-valued $\mathfrak{F}$-adapted L\'{e}vy process starting from $x\in \mathbb{R}^d$,
		with characteristics $(a,Q,\Pi)$.
		We introduce the convolution semigroup
		\begin{equation}
			\label{LevyConv}
			T_t\, f(x) = \mathbb{E} f(A^x_t) = \int_{\mathbb{R}^d} f(y) \mathbb{P}(A^x_t \in \mathrm dy)
		\end{equation}
		with infinitesimal generator
		\begin{equation}
			\label{adlera}
			\lim_{t \to 0} \frac{1}{t} \left( T_t\,f -f \right) = \mathcal{A}f,
		\end{equation} 
		where the strong limit exists in the domain
		\begin{equation}
			\label{cucitrice}
			D(\mathcal{A})  \mathrel{\mathop:}= \left\lbrace f \in L^1_{\text{loc}}(\mathbb{R}^d)
			\, :\, \int_{\mathbb{R}^d} |\hat{f}(\xi)|^2 \left( 1 + |\Psi(\xi)|^2 \right) \mathrm d\xi
			< \infty \right\rbrace.
		\end{equation}
		In \eqref{cucitrice}, $\hat{f}$ represents the Fourier transform of $f$.
		The cumulant generating function (also known in some literature as Fourier symbol) of the process $A_t := A_t^0$, $t \ge 0$, clearly is
		\begin{equation}
			\label{symbol-Levy}
			\Psi(\xi) = i \langle a, \xi \rangle + \frac{1}{2} \langle \xi, Q \xi \rangle +
			\int_{\mathbb{R}^d\setminus \{0\}} \left( 1- e^{i\langle z, \xi \rangle}
			+ i\langle z, \xi \rangle \mathbb{I}_{|z|<1} \right) \Pi( \mathrm dz).
		\end{equation}
		The Borel measure $\Pi(\cdot)$ is the so-called L\'evy measure satisfying
		$ \int_{\mathbb{R}^d\setminus \{0\}} (1 \wedge |z|^2) \, \Pi(dz)  < \infty$,
		where, as usual, $|z|^2 = \langle z, z \rangle$. We have that
		\begin{align}
			\label{Psi-multiplier}
			\mathcal{A}f (x) = \lim_{t \to 0} \frac{T_t \, f(x)- f(x)}{t}
			= \int_{\mathbb{R}^d} e^{i \langle x, \xi\rangle } \lim_{t \to 0} \frac{e^{-t \Psi(\xi)}
			- 1}{t} \hat{f}(\xi) \, \mathrm d\xi
			= \int_{\mathbb{R}^d} e^{i \langle x, \xi\rangle } \left( - \Psi(\xi) \right)
			\hat{f}(\xi) \, \mathrm d\xi
		\end{align}
		and therefore, $-\Psi$ is the Fourier multiplier of $\mathcal{A}$. An explicit representation in terms of L\'evy measure
		can be also given, we consider some special cases below.
		The semigroup $\mu_t(y, x) = \mathbb{P}(A^x_t \in \mathrm dy)/ \mathrm dy$ denotes the density of the
		L\'{e}vy process $A^x_t$ on $\mathbb{R}^d$ starting from $x \in \mathbb{R}^d$. This means that
		\begin{equation}
			\label{charLevy}
			\mathbb{E} \exp\left(i \xi A_t\right) = \int_{\mathbb{R}^d} e^{i\xi y} \mu_t(\mathrm dy)
			= \exp \left(-t \Psi(\xi) \right)
		\end{equation}
		and that the function $\Psi(\cdot)$ completely determines the density of $A_t$, $t \ge0$.
		In \eqref{charLevy}, $\mu_t(\mathrm dy) = \mu_t(\mathrm dy, 0)$.
		
		We introduce the time-change operator
		\begin{equation}
			H^{\gamma, \nu, \delta}_t = \int_0^\infty h(\mathrm dy, t) T_y \label{conv-TC-telegraph}
		\end{equation}
		where $h(y,t)$ is the solution to \eqref{04sub} and hence the density law of
		$Z_t^\delta$, $t \ge 0$, with $\delta \in (0, \infty)$
		and $T_y$ is the convolution semigroup in \eqref{LevyConv}. For
		$\delta \in (0,1)$, the function $h(y,t)$ coincides with the density law of		
		\eqref{Eproc}.
		
		For the operator \eqref{conv-TC-telegraph}, we obtain that
		\begin{equation}
			\| H^{\gamma, \nu, \delta}_t  \, f\| \leq  \int_0^\infty | h(\mathrm dy, t)| \, \| T_y\, f\|
			\leq  \| f\| \int_0^\infty | h(\mathrm dy, t)| = \| f\|. 
		\end{equation}
		Indeed, $T_t$ is a strongly continuous contraction semigroup on $\textup{C}^\infty(\mathbb{R}^d)$ and $h(\cdot, t)$ is a
		probability measure on $(0, \infty)$.
		Consider now the space $L^p((0,\infty), e^{-t} \mathrm dt)$ of all measurable functions $t \mapsto f(t)$
		equipped with the norm $\| f\|_p^p = \int_0^\infty |f(t)|^p e^{-t} \mathrm dt$.
        Recalling that $n=\lceil \delta \rceil$ and therefore that $\delta/n \in (0,1]$, we have
		\begin{equation}
			\| H^{\gamma, \nu, \delta}_t  \, f\|_1 \leq \int_0^\infty \| h(\mathrm dy, t)\|_1 \, |T_y \,
			f| =  \| T_{t/2^{\delta/n}} \, f\|_1.
		\end{equation} 
		
		Indeed, from \eqref{tetration} we also have that
		\begin{align}
			\| h(\mathrm dy, \cdot)\|_1 = \int_0^\infty h(\mathrm dy, t) e^{-t} \mathrm dt
			= 2^{\delta/n} e^{-y 2^{\delta/n} } \mathrm dy.
		\end{align}
		From Proposition \ref{mariolino} we have that $h(x, t) = \int_0^\infty l_{\delta/n}(x, r) k(\mathrm dr,t) $
		and the operator \eqref{conv-TC-telegraph} takes the form
		\begin{equation}
			\label{l-T-operator}
			H^{\gamma, \nu, \delta}_t = \int_0^\infty k(\mathrm dr, t)
			\int_0^\infty  l_{\delta/n}(\mathrm dy, r)\,T_y
		\end{equation}
		where, $l_{\delta/n}(\mathrm dy, r) =\mathbb{P}\{L^{\delta/n}_r \in \mathrm dy\}$
		and $k(\mathrm dr, t) = \mathbb{P} \{ F_t \in \mathrm dr \}$.

		For $n=1$ the density law of the process $F_t \overset{\text{d}}{=}
		C^{(\gamma+\nu, \nu)}_t$, $t \ge 0$ has the explicit representation (see \cite{DTO2012})
		\begin{equation}
			\label{explicit-h}
			k(x,t) = \int_0^t l_{\gamma+\nu}(x,y) v_{\nu}(t-y, x) \mathrm dy
			+ \int_0^t l_{\nu}(x, y) v_{\gamma+\nu}(t-y,x) \mathrm dy.
		\end{equation}
		For the sake of completeness we show that the Laplace transform of the density \eqref{explicit-h} is written as
		\begin{equation}
			\tilde{k}(x,s) = \frac{1}{s} \left(  s^{\gamma + \nu } + s^{\nu} \right) e^{-x(s^{\gamma+\nu} + s^\nu)}.
		\end{equation}
		Indeed, from the Laplace transforms $\tilde{l}_\alpha(x,s)=s^{\alpha-1}e^{-xs^\alpha}$
		and $\tilde{v}_\alpha(s, x) = e^{-xs^\alpha}$, we arrive at
		\begin{align}
			\int_0^\infty e^{-st}k(x,t)dt = \tilde{l}_{\gamma+\nu}(x, s)
			\tilde{v}_\nu(s, x) + \tilde{l}_\nu(x, s) \tilde{v}_{\gamma+\nu}(s, x) 
			= \left( s^{\gamma + \nu - 1} + s^{\nu -1} \right) e^{-x(s^{\gamma+\nu} + s^\nu)}.
		\end{align}
		Note also that when $\delta=n=0$, the process
		$F_t \overset{\text{d}}{=} L_t^{\gamma+\nu}$, which is
		an inverse $(\gamma+\nu)$-stable subordinator.
			
		Recall that $\mathbb{R}^{d+1} := \mathbb{R}^d\times(0,\infty)$ and that the function space 
        $\mathbb{A}^{k}$ is defined in \eqref{rafaela}.
   
		\begin{theorem}
			\label{TC-telegraph}
			Let $\delta > 0$, $n= \lceil \delta \rceil$, and $\delta \nu < \gamma+\nu \le 1$, $\gamma,\nu \in (0,1)$. The unique solution to 
			\begin{align}
				\label{Problem-L}
				\begin{cases}
                    g \in \mathbb{A}^{n}(\mathbb{R}^{d+1}),\\
					\left( \mathbb{D}_{\nu,\gamma+\nu,-1;0+}^\delta g(x,\cdot) \right) (t)
					= \mathcal{A}g(x, t), & x \in \mathbb{R}^d, \: t>0, \\
					g(x, 0) = f(x),
				\end{cases}
			\end{align}
            with $f \in D(\mathcal{A})$, is written as $g(x, t) = H^{\gamma, \nu, \delta}_t\, f(x) = \mathbb{E}f \Bigl(A^x_{Z_t^\delta} \Bigr)$,
			where $H^{\gamma, \nu, \delta}_t$ is the time-change operator \eqref{l-T-operator} and
			$A_t^x$, $t \ge 0$, is the L\'evy process started at $x \in \mathbb{R}^d$ with infinitesimal
			generator \eqref{adlera}.  

			\begin{proof}
                From \eqref{05lap} and \eqref{conv-TC-telegraph}
                we obtain the Laplace transform 
                \begin{align}
                    \tilde{g}(x, s) = & \int_0^\infty \tilde{h}
                    (\mathrm dy, s) T_y\, f(x)
                    = \int_0^\infty \tilde{h}(\mathrm dy, s)
                    \, \mathbb{E} f(A^x_y).
                \end{align}
                We need now the Fourier transform
                \begin{align}
                    \label{finestra}
                    \int_{\mathbb{R}^d} e^{i \xi \cdot x} \,
                    \tilde{g}(x, s)\, \mathrm dx = & \hat{f}(\xi)
                    \int_0^\infty \tilde{h}(\mathrm dy, s)
                    \, \hat{\mu}_y(\xi)
                \end{align}
                where recall that $\hat{\mu}_y(\xi)= e^{- y\Psi(\xi)}$.
                From \eqref{05lap} and \eqref{charLevy},
                \eqref{finestra} takes the form
                \begin{align}
                    \label{sam}
                    \hat{\tilde{g}}(\xi, s) & = \hat{f}(\xi)
                    \frac{s^{\gamma + \nu -1}
                    (1+s^{-\nu})^\delta}{\Psi(\xi)
                    + s^{\gamma + \nu }(1+s^{-\nu})^\delta}.
                \end{align}
                The Fourier transform of \eqref{Problem-L} leads to the equation
                \begin{equation}
                    \left( \mathbb{D}_{\nu,\gamma+\nu,-1;0+}^\delta \hat{g}(\xi,\cdot) \right) (t) = -\Psi(\xi)\, \hat{g}(\xi,t).
                \end{equation}
                By taking into account formula \eqref{LapCaputo-typeOp} we get the Fourier-Laplace transform
                \begin{align}
                    s^{\gamma+\nu}(1+s^{-\nu})^{\delta}\hat{\tilde{g}}(\xi, s)
                    - \hat{g}(\xi, 0^+) s^{\gamma+\nu-1}(1+s^{-\nu})^{\delta} + \Psi(\xi)\, \hat{\tilde{g}}(\xi,s)=0
                \end{align}
                and therefore, we get that
                \begin{align}
                    \hat{\tilde{g}}(\xi,s) = \hat{g}(\xi, 0^+) \frac{s^{\gamma + \nu
                    -1}(1+s^{-\nu})^\delta}{\Psi(\xi) + s^{\gamma + \nu }(1+s^{-\nu})^\delta}.
                \end{align}
                with $g(x,0)=f(x)$.
            \end{proof}
        \end{theorem}

        \begin{remark}
            Consider the integral \eqref{l-T-operator}. We observe that 
            \begin{equation}
                \int_0^\infty  l_{\delta/n}(\mathrm dy, r)\,T_y\, f(x) = \mathbb{E}f(A^x_{L^{\delta/n}_r})
            \end{equation}
            is the solution to the fractional problem
            \begin{align}
                \begin{cases}
                    \frac{\partial^{\delta/n}}{\partial r^{\delta/n}}u(x,r)
                    = \mathcal{A}u(x,r), \quad x\in \mathbb{R}^d,\; r>0, \; \delta/n \in (0,1),\\
                    u(x,0)=f(x).
                 \end{cases}
            \end{align} 
            Furthermore, for $\delta=1$, $\int_0^\infty k(\mathrm dr, t)T_r\, f(x) = \mathbb{E}f(A^x_t)$ is the solution to
            \begin{align}
                \begin{cases}
                    \frac{\partial^{\gamma+\nu}}{\partial t^{\gamma+ \nu}}u(x,r) +
                    \frac{\partial^\nu}{\partial t^\nu}u(x,r) = \mathcal{A}u(x,t), \quad x\in \mathbb{R}^d,
                    \; t>0, \; \nu < \gamma+\nu \le 1,\\
                u(x,0)=f(x).
                \end{cases}
            \end{align} 
        \end{remark}

        \begin{remark}
            Note that even though Theorem \ref{TC-telegraph} requires that $n \nu < \gamma+\nu \le 1$, the Fourier-Laplace
            transform \eqref{sam} is still valid for $\delta \nu < \gamma+\nu \le 2$. In this case however it cannot be
            related to the process $A^x_{Z_t^\delta}$, $t \ge 0$.
        \end{remark}

        \begin{remark}
            Formula \eqref{symbol-Levy} defines the Fourier multiplier of $\mathcal{A}$ which is the infinitesimal generator of
            a L\'{e}vy process $A_t$, $t \ge 0$.
            We mention below some specific cases:
            \begin{itemize}
                \item if $\Psi(\xi)=|\xi|^{2\alpha}$ with $\alpha \in (0,1]$, then $\mathcal{A} = -(-\triangle)^{\alpha}$ is the fractional Laplacian.
                    The process $A_t$ is an isotropic stable process and becomes a Brownian motion for $\alpha=1$.
                    Thus, for $\alpha \in (0,1)$, we have that (for a well-defined function $f$)
                    \begin{equation}
                        \label{frac-laplacian}
                        -\mathcal{A} f(x) = (-\triangle)^{\alpha} f(x) = \mathcal{C}(\alpha, d)
                        \int_{\mathbb{R}^d} \frac{f(x+y) + f(x-y)-2f(x)}{|y|^{2\alpha +d}} \mathrm dy,
                    \end{equation}
                    where $\mathcal{C}(\alpha, d)$ is a constant depending on $\alpha,d$.
                    It is well known that, in this case, the process $A_t = B_{V^\alpha_t}$ where $B$
                    is a multidimensional Brownian motion and $V^\alpha$ is a stable subordinator (the Bochner subordination rule).
                    In our case, therefore we get that $g(x,t) = \mathbb{E}f(x+B_{\tau_t})$,
                    where $\tau_t = V^\alpha_{Z^\delta_t}$ is a time-changed subordinator,
                    is the solution to the problem \eqref{Problem-L} with generator
                    \eqref{frac-laplacian}.
                \item if $d=1$ and $\Psi(\xi) = \lambda (1-e^{i\xi})$, then $A_t$ is a Poisson process on $\mathbb{Z}_{+}:=\{0,1,2, \ldots \}$
                    with rate $\lambda >0$ and 
                    \begin{equation}
                        \mathcal{A}f(x) = \lambda\{ f(x) - f(x-1)\} 1_{\mathbb{Z}_+}(x) 
                    \end{equation}
                    is the governing operator written as $\lambda$ times the discrete gradient on $\mathbb{Z}_+$.
                    Also, we usually write $\mathcal{A}f = \lambda(I-B)f$.
                \item if $d=1$ and $\Psi(\xi) = \lambda (1+i\xi - e^{i\xi})$, the corresponding process is the compensated Poisson on $\mathbb{R}$
                    with rate $\lambda >0$. The generator takes the form
                    \begin{equation}
                        \mathcal{A}f = \lambda (I-B)f - \lambda f^\prime.
                    \end{equation}
            \end{itemize}
        \end{remark}
	
	\section{One-dimensional case with $\bm{0 <\gamma+\nu \le 2}$}
	
        \label{last}
		Here the results obtained in the previous section are specialized for
		$\mathcal{A} = c \partial^2/\partial x^2$, $0 \ne c \in \mathbb{R}$.
        Note however that in this section the order of the operator
        $\mathbb{D}_{\nu,\gamma+\nu,-\lambda;0+}^\delta$
        is allowed to reach 2. Specifically we need $0< \gamma+\nu \le 2$.      
		
		Consider thus the Cauchy problem
		\begin{align}
			\label{00garlach}
			\begin{cases}
				\left( \mathbb{D}_{\nu,\gamma+\nu,-\lambda;0+}^\delta g(x,\cdot) \right) (t)
				= c \frac{\partial^2}{\partial x^2} g(x,t), \qquad x \in \mathbb{R}, \: t > 0, \\
				g(x,0^+) = \delta(x), \\
				\frac{\partial}{\partial t} g(x,t) \bigr|_{t \downarrow 0} = 0.
			\end{cases}
		\end{align}
		In the above equation
		$\delta \in \mathbb{R}$, $\gamma \in (0,\infty)$, $\nu \in (0,\infty)$,
		$0 < \gamma + \nu \le 2$, $c \ne 0$ is a real constant. Note therefore that here $\gamma+\nu \in (0,2]$
		so that \eqref{00garlach} is not in fact a direct specialization of \eqref{Problem-L}.
		This explains also the presence in \eqref{00garlach} of the addictional initial condition
		$\frac{\partial}{\partial t} g(x,t) \bigr|_{t \downarrow 0} = 0$. 
		
		\begin{remark}
			For $\delta = 0$, equation \eqref{00garlach} reduces to the time-fractional diffusion-wave equation
			\citep{ors1}.
			Indeed we have
			\begin{align}
				& \hspace{-1cm} \left( \mathbb{D}_{\nu,\gamma+\nu,-\lambda;0+}^0 g(x,\cdot) \right) (t)
				= c \frac{\partial^2}{\partial x^2} g(x,t) \\
				\Leftrightarrow {} & \quad
				\frac{\partial^{\gamma+\nu+\theta}}{\partial t^{\gamma+\nu+\theta}}
				\int_0^t (t-y)^{\theta-1} E_{\nu,\theta}^0 \left[ -\lambda (t-y)^\nu \right]
				g(x,y) \, \mathrm dy \notag
				= c \frac{\partial^2}{\partial x^2} g(x,t) + \delta(x) \, t^{-(\gamma+\nu)}
				E_{\nu,1-(\gamma+\nu)}^0 (-\lambda t^\nu) \notag \\
				\Leftrightarrow {} & \quad
				\frac{\partial^{\gamma+\nu}}{\partial t^{\gamma+\nu}}
				\frac{\partial^\theta}{\partial t^\theta} \frac{1}{\Gamma(\theta)} \int_0^t (t-y)^{\theta-1}
				g(x,y) \, \mathrm dy
				= c \frac{\partial^2}{\partial x^2} g(x,t) + \delta(x)
				\frac{ t^{-(\gamma+\nu)}}{\Gamma(1-(\gamma+\nu))} \notag \\
				& \hspace{-.9cm} \overset{(\kappa=\gamma+\nu)}{\Leftrightarrow} \frac{\partial^\kappa}{\partial t^\kappa}
				g(x,t) = c \frac{\partial^2}{\partial x^2} g(x,t) + \delta(x)
				\frac{ t^{-\kappa}}{\Gamma(1-\kappa)}, \qquad 0<\kappa \le 2 \notag \\
				\Leftrightarrow {} & \quad \frac{\mathfrak{d}^\kappa}{\mathfrak{d} t^\kappa}
				g(x,t) = c \frac{\partial^2}{\partial x^2} g(x,t). \notag
			\end{align}
			In the second step of the above formula we have used the fact that
			\begin{align}
				E_{\nu,\theta}^0 \left[ -\lambda (t-y)^\nu \right] = 1/\Gamma(\theta), \qquad
				E_{\nu,1-(\gamma+\nu)}^0 (-\lambda t^\nu) = 1/\Gamma(1-(\gamma+\nu)).
			\end{align}
			Also, as mentioned before, we considered here that the semigroup property for the Riemann--Liouville
			fractional derivative and hence some regularity conditions on the solution $g$ are fulfilled
			(see Section \ref{frank}).
		\end{remark}
		
		\begin{remark}
			For $\delta = 1$, equation \eqref{00garlach} reduces to the fractional telegraph equation
			\citep{ors2}.
			In this case we have
			\begin{align}
				& \hspace{-1cm}
				\left( \mathbb{D}_{\nu,\gamma+\nu,-\lambda;0+}^1 g(x,\cdot) \right) (t)
				= c \frac{\partial^2}{\partial x^2} g(x,t) \\
				\Leftrightarrow {} & \quad
				\frac{\partial^{\gamma+\nu+\theta}}{\partial t^{\gamma+\nu+\theta}}
				\int_0^t (t-y)^{\theta-1} E_{\nu,\theta}^{-1} \left[ -\lambda (t-y)^\nu \right]
				g(x,y) \, \mathrm dy
				= c \frac{\partial^2}{\partial x^2} g(x,t) + \delta(x) \, t^{-(\gamma+\nu)}
				E_{\nu,1-(\gamma+\nu)}^{-1} (-\lambda t^\nu). \notag
			\end{align}
			Now, by considering that
			\begin{align}
				& E_{\nu,\theta}^{-1} \left[ -\lambda (t-y)^\nu \right]
				= \frac{1}{\Gamma(\theta)} + \frac{\lambda (t-y)^\nu}{\Gamma(\nu+\theta)}, \\
				& E_{\nu,1-(\gamma+\nu)}^{-1} (-\lambda t^\nu)
				= \frac{1}{\Gamma(1-(\gamma+\nu))} + \frac{\lambda t^\nu}{\Gamma(1-\gamma)},
			\end{align}
			we can write
			\begin{align}
				\frac{\partial^{\gamma+\nu}}{\partial t^{\gamma+\nu}} &
				\frac{\partial^\theta}{\partial t^\theta} \frac{1}{\Gamma(\theta)} \int_0^t (t-y)^{\theta-1}
				g(x,y) \, \mathrm dy
				+\lambda  \frac{\partial^\gamma}{\partial t^\gamma} \frac{\partial^{\nu+\theta}}{\partial^{\nu+\theta}}
				\frac{1}{\Gamma(\nu+\theta)} \int_0^t (t-y)^{\nu+\theta-1} g(x,y) \, \mathrm dy \\
				& = c \frac{\partial^2}{\partial x^2} g(x,t) + \delta(x)
				\frac{ t^{-(\gamma+\nu)}}{\Gamma(1-(\gamma+\nu))} + \delta(x)
				\frac{\lambda t^{-\gamma}}{\Gamma(1-\gamma)} \notag \\
				\Leftrightarrow \quad {} & \frac{\partial^{\gamma+\nu}}{\partial t^{\gamma+\nu}}
				g(x,t) + \lambda \frac{\partial^\gamma}{\partial t^\gamma} g(x,t)
				= c \frac{\partial^2}{\partial x^2} g(x,t) + \delta(x)
				\frac{ t^{-(\gamma+\nu)}}{\Gamma(1-(\gamma+\nu))} + \delta(x)
				\frac{\lambda t^{-\gamma}}{\Gamma(1-\gamma)} \notag \\
				\Leftrightarrow \quad {} & \frac{\mathfrak{d}^{\gamma+\nu}}{\mathfrak{d} t^{\gamma+\nu}}
				g(x,t) + \lambda \frac{\mathfrak{d}^\gamma}{\mathfrak{d} t^\gamma} g(x,t)
				= c \frac{\partial^2}{\partial x^2} g(x,t), \notag
			\end{align}
			where $0< \gamma+\nu \le 2$.
		\end{remark}
		
		The Laplace--Fourier transform of the solution $g(x,t)$ to equation \eqref{00garlach} can be easily determined
		as follows. We start with the application of the Fourier transform
		$\hat{g}(\beta, t) = \int_{-\infty}^\infty e^{i\beta x} g(x,t) \, \mathrm dx$,
		immediately yielding
		\begin{align}
			\label{20mm}
			\begin{cases}
				\left( \mathbb{D}_{\nu,\gamma+\nu,-\lambda;0+}^\delta \hat{g}(\beta,\cdot) \right) (t)
				= -c \beta^2 \hat{g}(\beta,t), \\
				\hat{g}(\beta,0^+) = 1, \\
				\frac{\partial}{\partial t} \hat{g}(x,t) \bigr|_{t \downarrow 0} = 0.
			\end{cases}
		\end{align}
		Then, by using formula \eqref{20de} and applying the Laplace transform (with parameter $s$) to both members of \eqref{20mm},
		the complete Laplace--Fourier transform
		of \eqref{00garlach} can be written as
		\begin{align}
			s^{\gamma+\nu+\theta} \hat{\tilde{g}}(\beta,s) \, s^{-\theta} \left( 1+\lambda s^{-\nu} \right)^\delta
			= -c \beta^2 \hat{\tilde{g}}(\beta,s) + s^{\gamma+\nu-1} \left( 1+\lambda s^{-\nu} \right)^\delta.
		\end{align}
		From this we immediately obtain the Laplace--Fourier transform of the solution to equation
		\eqref{00garlach} as
		\begin{align}
			\label{01lf}
			\hat{\tilde{g}}(\beta,s) = \frac{s^{\gamma+\nu-1}
			\left( 1+\lambda s^{-\nu} \right)^\delta}{s^{\gamma+\nu} \left( 1+\lambda s^{-\nu} \right)^\delta
			+c \beta^2}.
		\end{align}
		
		\begin{remark}
			Clearly, for $\delta=0$, formula \eqref{01lf} reduces to the Laplace--Fourier transform of the solution
			to the fractional diffusion-wave equation \citep[formula (2.17)]{ors1},
			while for $\delta=1$ it leads to that of the
			fractional telegraph equation \citep[formula (2.6) for $\gamma=\nu=\alpha$]{ors2}.
		\end{remark}
		
		The Fourier transform of the solution to \eqref{00garlach} can be derived by inverting the
		Laplace transform in \eqref{01lf} as follows.
		\begin{align}
			\hat{\tilde{g}}(\beta,s) = s^{-1} \left( 1+\frac{c\beta^2}{s^{\gamma+\nu}
			\left( 1+\lambda s^{-\nu} \right)^\delta} \right)^{-1}
			= s^{-1} \sum_{r=0}^\infty \left[ - \frac{c\beta^2}{s^{\gamma+\nu}
			\left( 1+\lambda s^{-\nu} \right)^\delta} \right]^r.
		\end{align}
		The last step is valid whenever $|c\beta^2/(s^{\gamma+\nu}(1+\lambda s^{-\nu})^\delta)|<1$.
		We then have
		\begin{align}
			\hat{\tilde{g}}(\beta,s) = \sum_{r=0}^\infty (-c\beta^2)^r s^{-(r(\gamma+\nu)+1)} \left(
			1+\lambda s^{-\nu}\right)^{-r\delta}.
		\end{align}
		Consequently, by recalling again formula \eqref{20de} and considering
		\citet[Theorem 30.1]{doetsch} which ensures the inversion term by term,
		the Fourier transform of $g(x,t)$ reads
		\begin{align}
			\label{02fourier}
			\hat{g}(\beta,t) & = \sum_{r=0}^\infty (-c\beta^2 t^{\gamma+\nu})^r
			E_{\nu,r(\gamma+\nu)+1}^{r\delta} \left( -\lambda t^\nu \right) \\
			& = \sum_{m=0}^\infty (-\lambda t^\nu)^m {}_2\psi_2 \left[ c\beta^2 t^{\gamma+\nu}
			\left|
			\begin{array}{l}
				(1,1),(m,\delta) \\
				(0,\delta),(\gamma+\nu,\nu m+1)
			\end{array}
			\right. \right] \notag,
		\end{align}
		where ${}_p \psi_q$ is the generalized Wright function \citep[Section 1.11]{kilbas}, defined as
		\begin{align}
			{}_p \psi_q (x) = {}_p \psi_q \left[ x \left|
			\begin{array}{l}
				(a_1,\alpha_1),\dots,(a_p,\alpha_p) \\
				(b_1,\beta_1),\dots,(b_q,\beta_q)
			\end{array}
			\right. \right]
			= \sum_{k=0}^\infty \frac{x^k}{k!} \frac{\prod_{m=1}^p \Gamma(a_m+\alpha_m k)}{\prod_{j=1}^q
			\Gamma(b_j+\beta_j k)},
		\end{align}
		where $x,a_m,b_j \in \mathbb{C}$, $\alpha_m,\beta_j \in \mathbb{R}$, $m=1,\dots, p$, $j=1,\dots,q$.
		
		The Laplace--Fourier transform \eqref{01lf} immediately yields
		\begin{align}
			\tilde{g}(x,s) & = \int_0^\infty e^{-st} g(x,t) \, \mathrm dt
			= \frac{1}{2\pi} \int_{-\infty}^\infty
			e^{-i\beta x} \frac{s^{\gamma+\nu-1} \left( 1+\lambda s^{-\nu} \right)^\delta}{
			s^{\gamma-\nu} \left( 1+\lambda s^{-\nu} \right)^\delta +c\beta^2} \, \mathrm d\beta \\
			& = \frac{1}{2sc^{1/2}} s^{(\gamma+\nu)/2} \left( 1+\lambda s^{-\nu} \right)^{\delta/2}
			e^{-\frac{|x|}{c^{1/2}}s^{(\gamma+\nu)/2} \left( 1+\lambda s^{-\nu} \right)^{\delta/2}}. \notag
		\end{align}
				
		\begin{remark}
			Note how the Fourier transform \eqref{02fourier}, for $\delta=0$ reduces to the
			Fourier transform of the solution of the pure fractional diffusion equation as we
			recall that $(0)_m = 0$, $m=1,2,\dots$, but $(0)_0 = 1$.
			Therefore,
			\begin{align}
				\hat{g}(\beta,t) = E_{\gamma+\nu,1}\left( -c\beta^2 t^{\gamma+\nu} \right),
			\end{align}
			which coincides with the corresponding formula of \citet{ors1}, page 212.
			
			When $\delta=1$ instead we arrive at
			\begin{align}
				\hat{g}(\beta,t) & = \sum_{r=0}^\infty (-c\beta^2 t^{\gamma+\nu})^r
				E_{\nu,r(\gamma+\nu)+1}^r \left( -\lambda t^\nu \right) \\
				& = \sum_{m=0}^\infty (-\lambda t^\nu)^m {}_2\psi_2 \left[ c\beta^2 t^{\gamma+\nu}
				\left|
				\begin{array}{l}
					(1,1),(m,1) \\
					(0,1),(\gamma+\nu,\nu m+1)
				\end{array}
				\right. \right] \notag,
			\end{align}
			which should be compared with formula (2.7) of \citet{ors2} when $\nu=\gamma=\alpha$.
		\end{remark}
		
		\begin{remark}
			For $\delta=2$, $\gamma>\nu$, we obtain an interesting specific case. In this case equation \eqref{00garlach}
			reduces to
			\begin{align}
				\begin{cases}
					\frac{\mathfrak{d}^{\gamma+\nu}}{\mathfrak{d} t^{\gamma+\nu}} g(x,t) + 2 \lambda
					\frac{\mathfrak{d}^\gamma}{\mathfrak{d} t^\gamma} g(x,t)
					+ \lambda^2 \frac{\mathfrak{d}^{\gamma-\nu}}{\mathfrak{d} t^{\gamma-\nu}}
					g(x,t)
					= c \frac{\partial^2}{\partial x^2} g(x,t), \\
					g(x,0^+) = \delta(x).
				\end{cases}
			\end{align}
			
			In the general case of $\delta = n \in \mathbb{N} \cup \{ 0 \}$, we can work out equation \eqref{00garlach}
			as follows.
			\begin{align}
				\left( \mathbb{D}_{\nu,\gamma+\nu,-\lambda;0+}^n g(x,\cdot) \right) (t)
				= c \frac{\partial^2}{\partial x^2} g(x,t)
				\Leftrightarrow {}  \quad &
				\frac{\partial^{\gamma+\nu+\theta}}{\partial t^{\gamma+\nu+\theta}}
				\int_0^t (t-y)^{\theta-1} E_{\nu,\theta}^{-n} \left[ -\lambda (t-y)^\nu \right]
				g(x,y) \, \mathrm dy  \\
				& = c \frac{\partial^2}{\partial x^2} g(x,t) + \delta(x) \, t^{-(\gamma+\nu)}
				E_{\nu,1-(\gamma+\nu)}^{-n} (-\lambda t^\nu). \notag
			\end{align}
			Now, considering that
			$(-n)_r = (-1)^r (n-r+1)_r = (-1)^r n!/(n-r)!$,
			we can write
			\begin{align}
				\label{22tt}
				& \frac{\partial^{\gamma+\nu+\theta}}{\partial t^{\gamma+\nu+\theta}}
				\int_0^t (t-y)^{\theta-1} \sum_{r=0}^n \binom{n}{r} \frac{\lambda^r(t-y)^{\nu r}}{\Gamma(\nu r+\theta)}
				g(x,y) \, \mathrm dy
				 = c \frac{\partial^2}{\partial x^2} g(x,t) + \delta(x)
				\sum_{r=0}^n \binom{n}{r} \frac{\lambda^r t^{-\left(\gamma-\nu\left(r-1\right)\right)}}{\Gamma \left(1-
				\left(\gamma-\nu\left(r-1\right)\right)\right)}  \\
				& \Leftrightarrow \quad
				\sum_{r=0}^n \binom{n}{r} \lambda^r \frac{\partial^{\gamma-\nu(r-1)}}{\partial
				t^{\gamma-\nu(r-1)}} g(x,t)
				= c \frac{\partial^2}{\partial x^2} g(x,t) + \delta(x)
				\sum_{r=0}^n \binom{n}{r} \frac{\lambda^r t^{-\left(\gamma-\nu\left(r-1\right)\right)}}{\Gamma \left(1-
				\left(\gamma-\nu\left(r-1\right)\right)\right)} \notag \\
				& \Leftrightarrow \quad
				\sum_{r=0}^n \binom{n}{r} \lambda^r \frac{\mathfrak{d}^{\gamma-\nu(r-1)}}{\mathfrak{d}
				t^{\gamma-\nu(r-1)}} g(x,t)
				= c \frac{\partial^2}{\partial x^2} g(x,t). \notag
			\end{align}
			Here we have $0 < \gamma-\nu(r-1) \le 2$ so that $n\nu < \gamma+\nu \le 2$.
		\end{remark}
		
		\begin{remark}[Wave-telegraph equation]
			When $\gamma=\nu=1$, the equation considered is a wave-telegraph equation.
			In this case, the allowed range for $\delta$ is $\delta \le 1$. The interpolating equation
			reads
			\begin{align}
				\label{30garlach}
				\left( \mathbb{D}_{1,2,-\lambda;0+}^\delta g(x,\cdot) \right) (t)
				= c \frac{\partial^2}{\partial x^2} g(x,t), \qquad x \in \mathbb{R}, \: t \ge 0.
			\end{align}
		\end{remark}

		\begin{appendices}
		
		\section{Properties of the operator $\bm{\mathbb D_{\alpha,\eta,\zeta;0+}^\xi}$}
			\label{frank}

			We analyze here some properties of the operator $\mathbb D_{\alpha,\eta,\zeta;0+}^\xi$.
			In the following proposition we show that
			the integral operator \eqref{20ciao} is in fact a generalization of the
			left sided Riemann--Liouville fractional integral, to which it reduces for $\xi=0$. The convolution
			kernel here is no longer a power-law but it is instead a generalized Mittag--Leffler kernel.
			It follows furthermore that the operator \eqref{00tex} represents a
			generalization to the left sided Riemann--Liouville
			fractional derivative.

			\begin{proposition}
				\label{55prop1}
				For $\xi=0$ we have that
				\begin{align}
					\label{first}
					& \left( \bm{\mathrm E}_{\alpha,\eta,\zeta;0+}^0 f(\cdot) \right) (t) =
					\frac{1}{\Gamma(\eta)} \int_0^t (t-y)^{\eta-1} f(y) \mathrm dy, \\
					\label{second}
					& \left( \bm{\mathrm D}_{\alpha,\eta,\zeta;0+}^0 f(\cdot) \right) (t) =
					\frac{\mathrm d^{\eta}}{\mathrm dt^{\eta}} f(t),
				\end{align}
				with $\eta > 0$, $\alpha>0$, $\zeta \in \mathbb{R}$, $t \ge 0$.
			
				\begin{proof}
					Formula \eqref{first} can be derived by simply considering that		
					$E_{\alpha,\eta}^0 \left[ \zeta (t-y)^\alpha \right] = 1/\Gamma(\eta)$.
					In order to briefly prove formula \eqref{second} we can write that
					\begin{align}
						\label{20tex}
						\left( \bm{\mathrm D}_{\alpha,\eta,\zeta;0+}^0 f(\cdot) \right) (t)
						& = \frac{\mathrm d^{\eta+\theta}}{\mathrm dt^{\eta+\theta}}
						\int_0^t (t-y)^{\theta-1} E_{\alpha,\theta}^0\left[ \zeta (t-y)^\alpha \right]
						f(y) \, \mathrm dy \\
						& = \frac{\mathrm d^{\eta}}{\mathrm dt^{\eta}}
						\frac{\mathrm d^{\theta}}{\mathrm dt^{\theta}}
						\frac{1}{\Gamma(\theta)} \int_0^t (t-y)^{\theta-1} f(y) \, \mathrm dy
						= \frac{\mathrm d^{\eta}}{\mathrm dt^{\eta}} f(t), \notag
					\end{align}
					as the Riemann--Liouville fractional derivative is the left-inverse operator to the Riemann--Liouville
					fractional integral. As before, in the second step of formula \eqref{20tex} we considered that
					$E_{\alpha,\theta}^0 \left[ \zeta (t-y)^\alpha \right] = 1/\Gamma(\theta)$,
					and that the semigroup property for the Riemann--Liouville fractional derivative is fulfilled as
					\begin{align}
						\frac{\mathrm d^r}{\mathrm dt^r} \int_0^t (t-y)^{\theta-1}
						f(y) \, \mathrm dy\,  \biggr|_{t=0} = 0, \qquad \forall r \in \mathbb{N} \cup \{ 0 \}.
					\end{align}
				\end{proof}
			\end{proposition}		
			
			\begin{remark}
				\label{valmora}
				From Proposition \eqref{55prop1} it immediately follows that the operator
				\eqref{vtvt} represents a generalization
				of the left sided Dzhrbashyan--Caputo fractional derivative
				$\mathfrak{d}^\eta/\mathfrak{d}t^\eta$. Indeed, for $\xi=0$, we can write
				\begin{align}
					\left( \mathbb{D}_{\alpha,\eta,\zeta;0+}^0 f(\cdot) \right) (t)
					& = \left( \bm{\mathrm D}_{\alpha,\eta,\zeta;0+}^0 f(\cdot) \right) (t)
					- f(0^+) \, t^{-\eta} E_{\alpha,1-\eta}^0 (\zeta t^\alpha), \\
					& = \frac{\mathrm d^{\eta}}{\mathrm dt^{\eta}} f(t) - f(0^+) \frac{t^{-\eta}}{\Gamma(1-\eta)}
					= \frac{\mathfrak{d}^{\eta}}{\mathfrak{d}t^{\eta}} f(t). \notag
				\end{align}
				For more information on Dzhrbashyan--Caputo derivatives, the reader can refer for example to \citet{kilbas}.
			\end{remark}
		
			\begin{proposition}		

				The operator $\bm{\mathrm D}_{\alpha,\eta,\zeta;0+}^\xi$,
				$\eta > 0$, $\alpha>0$, $\zeta \in \mathbb{R}$, $t \ge 0$,
				can be written also as
				\begin{align}
					\left( \bm{\mathrm D}_{\alpha,\eta,\zeta;0+}^\xi f( \cdot ) \right) (t) =
					W_{-1,\xi+1} \left( -\zeta J^\alpha_t \right) f(t),
				\end{align}
				where
				\begin{align}
					W_{a,b} (x)
					= \sum_{r=0}^\infty \frac{x^r}{r!\Gamma(a r+b)}, \qquad a \ge 1, b \in \mathbb{R},
				\end{align}
				is the classical Wright function which is convergent in $|x|<1$ if $a=-1$, $b>0$,
				and $ J^\alpha_t$ is the Riemann--Liouville fractional integral.
		
				\begin{proof}
		
					We start by expanding in series the generalized Mittag--Leffler function in the kernel
					of the operator.
					\begin{align}
						\left( \bm{\mathrm D}_{\alpha,\eta,\zeta;0+}^\xi f( \cdot ) \right) (t)
						& = \frac{\mathrm d^{\eta+\theta}}{\mathrm d t^{\eta+\theta}}
						\int_0^t (t-y)^{\theta-1} E_{\alpha,\theta}^{-\xi}\left[ \zeta (t-y)^\alpha \right]
						f(y) \, \mathrm dy \\
						& = \frac{\mathrm d^{\eta+\theta}}{\mathrm d t^{\eta+\theta}} \int_0^t (t-y)^{\theta-1}
						\sum_{r=0}^\infty \frac{\zeta^r (-\xi)_r (t-y)^{\alpha r}}{r! \Gamma(\alpha r + \theta)}
						f(y) \, \mathrm dy. \notag
					\end{align}
					By recalling now that $(-\xi)_r = (-1)^r (\xi-r+1)_r = (-1)^r \Gamma(\xi+1)/\Gamma(\xi-r+1)$
					and considering again that
					the semigroup property is satisfied,
					we have that
					\begin{align}
						\left( \bm{\mathrm D}_{\alpha,\eta,\zeta;0+}^\xi f(\cdot) \right)(t)
						& = \Gamma(\xi+1) \frac{\mathrm d^{\eta+\theta}}{\mathrm dt^{\eta+\theta}}
						\sum_{r=0}^\infty \frac{(-\zeta)^r}{r!\Gamma(\xi-r+1)}
						\frac{1}{\Gamma(\alpha r+\theta)} \int_0^t (t-y)^{\alpha r+\theta -1} f(y) \, \mathrm dy \\
						& = \Gamma(\xi+1) \frac{\mathrm d^{\eta+\theta}}{\mathrm d t^{\eta+\theta}}
						\sum_{r=0}^\infty \frac{(-\zeta)^r}{r!\Gamma(\xi-r+1)} J^{\alpha r+\theta}_t f(t) \notag \\
						& = \Gamma(\xi+1) \frac{\mathrm d^{\eta}}{\mathrm d t^{\eta}}
						\sum_{r=0}^\infty \frac{(-\zeta)^r}{r!\Gamma(\xi-r+1)} J^{\alpha r}_t f(t), \notag
					\end{align}
					provided that the series converges and
					where $J^{\alpha r}_t$ represents the Riemann--Liouville fractional integral of order $\alpha r$.
					Thus we obtain that
					\begin{align}
						\label{perepe}
						\left( \bm{\mathrm D}_{\alpha,\eta,\zeta;0+}^\xi f(\cdot) \right)(t)
						& = \Gamma(\xi+1) \frac{\mathrm d^{\eta}}{\mathrm d t^{\eta}}
						\sum_{r=0}^\infty \frac{(-\zeta J^\alpha_t)^r}{r!\Gamma(\xi-r+1)} f(t) \\
						& = \Gamma(\xi+1) \frac{\mathrm d^{\eta}}{\mathrm d t^{\eta}}
						W_{-1,\xi+1} \left( -\zeta J^\alpha_t \right) f(t). \notag
					\end{align}
					The obtained representation \eqref{perepe} is formal and it becomes an actual representation whenever
					all the requested convergence conditions are fulfilled.
				\end{proof}	
			\end{proposition}
			
		\end{appendices}

	\bibliography{mirko} 
	\bibliographystyle{plainnat}

\end{document}